\documentclass[12pt]{amsart}

\usepackage{amssymb,amsfonts,amsthm,amsmath}

\usepackage{latexsym,multicol,graphicx}

\usepackage{graphicx, verbatim}

\usepackage{xcolor}
\usepackage{verbatim}

\usepackage[a4paper, hmargin=2.5cm, vmargin={3.5cm, 3.5cm}]{geometry}

\usepackage{fancyhdr}

\newtheorem*{theorem*}{Theorem}
\newtheorem{lemma}{Lemma}

\renewcommand{\phi}{\varphi}

\renewcommand{\emptyset}{\varnothing}

\newcommand{\e}{\varepsilon}

\newcommand{\eps}{\varepsilon}

\newcommand{\la}{\lambda}

\newcommand{\cH}{\mathcal H}

\newcommand{\cP}{\mathbb P}

\newcommand{\bE}{\mathbb E}

\newcommand{\vol}{\operatorname{vol}}

\newcommand{\sgn}{\operatorname{sgn}}

\renewcommand{\le}{\leqslant}
\renewcommand{\ge}{\geqslant}

\newcommand{\bR}{\mathbb R}
\newcommand{\bC}{\mathbb C}
\newcommand{\bZ}{\mathbb Z}
\newcommand{\bD}{\mathbb D}

\newcommand{\bP}{\mathbb P}
\newcommand{\bS}{\mathbb S}

\newcommand{\Cr}{\operatorname{Cr}}
\newcommand{\Var}{\operatorname{Var}}

\numberwithin{equation}{subsection}

\author{Fedor Nazarov}
\address{Department of Mathematics, Kent State University, Kent OH 44242, USA
}
\email{nazarov@math.kent.edu}

\author{Mikhail Sodin}
\address{School of Mathematical Sciences\\
Tel Aviv University\\
Tel Aviv 69978\\
Israel}
\email{sodin@tauex.tau.ac.il}

\thanks{This work was partially supported
by U.S. NSF Grant DMS-1900008 (F.N),
and by ERC Advanced Grant 692616
(M.S.)}

\begin{document}

{\mbox{}\hfill \today}

\bigskip

\title{Fluctuations in the number of nodal domains}

\begin{abstract}

We show that the variance of the number of connected components of the
zero set of the two-dimensional Gaussian ensemble of random spherical
harmonics of degree $n$ grows as a positive power of $n$.
The proof uses no special properties of spherical harmonics and works
for any sufficiently regular ensemble of Gaussian
random functions on the two-dimensional sphere with distribution invariant with respect
to isometries of the sphere.

Our argument connects the fluctuations in the number of nodal lines with those in
a random loop ensemble on planar graphs of degree four, which can be viewed as a step towards justification of the Bogomolny-Schmit heuristics.
\end{abstract}

\maketitle

\mbox{}\hfill{\sf In memory of Jean Bourgain}

\section{Introduction}

Let $(f_n)$ be the ensemble of random Gaussian spherical harmonics
of degree $n$ on the two-dimensional sphere, and let $N(f_n)$ be the number
of connected components of the zero set $\{ f_n=0\}$. It is known that
\[
\bE [N(f_n)] = (c+o(1))n^2,  \qquad n\to\infty,
\]
with a positive numerical constant $c$ and that the random variable $N(f_n)$
exponentially concentrates around its mean~\cite{NS1}.
A beautiful Bogomolny-Schmit heuristics~\cite{BS} suggests that, for any $\eps>0$ and
$n$ large enough,
\[
n^{2-\eps} < \Var[N(f_n)] < n^{2+\eps}.
\]
However, the rigorous bounds we are aware of are much weaker:
\[
%\begin{equation}\label{eq:1-10}
n^\sigma \lesssim \Var[N(f_n)] \lesssim n^{4-\sigma}
%\end{equation}
\]
with some $\sigma>0$. The upper bound with $\sigma=\frac2{15}$
follows from the exponential concentration of $N(f_n)$ around its
mean (see~\cite[Remark~1.2]{NS1}). The purpose of this paper is to prove
the lower bound. The proof we give uses no special
properties of spherical harmonics and shows that this lower
bound holds for any smooth, non-degenerate ensemble of Gaussian random functions
on the two-dimensional sphere $\bS^2$
with distribution invariant with respect to isometries of the sphere and
correlations decaying at least as a positive power of the appropriately
scaled distance on $\bS^2$.

It is worth mentioning that in a recent work~\cite{BMM} Beliaev, McAuley, Muirhead
found non-trivial
lower bounds for fluctuations of the number of connected components in the disk of radius
$R\gg 1$ of the level sets $\{F=\ell\}$
and of the excursion set $\{F\ge \ell\}$  of the random plane wave (which is a scaling limit of the
ensemble of spherical harmonics) for non-zero levels $\ell\ne 0$. It is expected that in their case the fluctuations are much larger than the ones we study. The techniques used in their work are quite different.

\section{The set-up and the main result}

Let $(f_L)$ be an ensemble of Gaussian random functions
on the two-dimensional sphere.
It is convenient to assume
that the function $f_L$ is defined on the sphere
$\bS^2(L)=\{x\in\bR^3\colon |x|=L\}$
of large radius $L$ and is normalized by $\bE[f_L(x)^2]=1$ for all $x\in \bS^2(L)$.
We always assume that the distribution of $f_L$ is invariant with respect to
the isometries of the sphere. Then the covariance kernel of $f_L$ has the form
\[
K_L(x, y) = \bE[f_L(x)f_L(y)] = k_L(d_L(x, y)), \qquad x, y\in\bS^2(L),
\]
where $d_L$ is the spherical distance on $\bS^2(L)$.
We call such an ensemble $(f_L)$ {\em regular} if the following two conditions hold:
\begin{enumerate}
\item{\em $\mathcal C^{3+}$-smoothness}: $K_L\in C^{3+\nu, 3+\nu}(\bS^2(L))$ with estimates uniform in $L$ and with some $\nu>0$.
\item{\em Power decay of correlations}: $K_L(x, y) \lesssim (1+d_L(x, y)\,)^{-\gamma}$, $x, y\in\bS^2(L)$, with some $\gamma>0$ and with the implicit constant independent of $L$.
\end{enumerate}
Condition (1) yields that almost surely $f_L\in C^3(\bS^2(L))$
with estimates uniform in $L$. We also note that condition (2)
is equivalent to the estimate $|k_L(d)|\lesssim (1+d\,)^{-\gamma}$ for $0\le d \le \pi L$ (with the implicit constant independent of $L$).

\medskip
By $Z(f_L)$ we denote the random zero set of $f_L$, which is, almost surely,  a collection of disjoint simple
smooth closed random curves (``loops'') on $\bS^2(L)$.
By $N(f_L)$ we denote the number of these loops.
\begin{theorem*}\label{thm:main}
Let $(f_L)$ be a regular Gaussian ensemble. Then there exists
$\sigma>0$ such that, for $L\ge L_0$,
\[
\Var[N(f_L)] \ge L^{\sigma}.
\]
\end{theorem*}

\medskip

There are many natural regular Gaussian ensembles, but a nuisance is that
the spherical harmonics ensemble is not among them. Spherical harmonics are symmetric with
respect to the center of the sphere so their values at the antipodal points on the sphere coincide up to the sign.
The correlations for this ensemble still satisfy condition (2) but only in the range
$0\le d \le (\pi-\e)L$ with any $\e>0$. For this reason, our theorem cannot be applied to this ensemble directly.
Luckily, the case of the spherical harmonics requires only
minor modifications in the proof of the theorem, which we will outline in the last section of this work. Essentially, we just need an analogue of our theorem for the projective plane $\bR\bP^2$ instead of the sphere.

\section{Main steps in the proof}

Heuristically, the fluctuations in the topology of the zero set are caused by fluctuations
in the signs of the critical values. To exploit this heuristics,
we fix a random function $f_L$
and slightly perturb it by a multiple of its independent copy $g_L$, i.e., consider the random function
\[
\widetilde{f}_L = \sqrt{1-\alpha'^2} f_L + \alpha' g_L,
\qquad 0<\alpha' \ll 1,
\]
which has the same distribution as $f_L$. Let $\alpha$ be another small parameter,
which is significantly bigger than $\alpha'$, $\alpha' \ll \alpha \ll 1$, and let
\[
\Cr(\alpha) = \bigl\{p\in \bS^2(L)\colon \nabla f_L(p)=0,\  |f_L(p)|\le \alpha \bigr\}.
\]
Then, as we will see, with high probability, given $f_L$,
the topology of the zero set
$Z(\widetilde{f}_L)$ is determined by the signs of $\widetilde{f}_L$ at $\Cr(\alpha)$,
that is, by the collection of the random values
$\bigl\{g_L(p)\colon p\in\Cr(\alpha) \bigr\}$. To make the correlations between these random values negligible, the set $\Cr(\alpha)$ should be well-separated on the sphere $\bS^2(L)$.
At the same time, the set $\Cr(\alpha)$ has to be relatively large; otherwise,
the impact of fluctuations in signs of $\widetilde{f}_L$ on the number
$ N(\widetilde{f}_L) $ will be negligible.
In Lemma~\ref{lemma_Cr(alpha)}, we will show that
\begin{itemize}
\item
{\em there exist positive $\eps_0$, $c$, and $C$ such that, given $0<\eps\le \eps_0$ and
$L\ge L_0(\e)$, for $L^{-2+\eps}\le\alpha\le L^{-2+2\eps}$, with probability very close to $1$, the set $\Cr(\alpha)$ is $L^{1-C\eps}$-separated and
$|\Cr(\alpha)|\ge L^{c\eps}$}.
\end{itemize}
The proof of this lemma given in Sections~7--9
is the longest and probably the most delicate part of our work.

\medskip
To understand how the signs of $\widetilde{f}_L$ at $\Cr(\alpha)$
affect the topology of the zero set $Z(\widetilde{f}_L)$, we develop in Section~6
a little caricature of the quantitative Morse theory. This caricature is non-random
- its applicability to the random function $f_L$ relies on the fact that with high
probability the Hessian $\nabla^2 f_L$ cannot degenerate at the points where the
function $f_L$ and its gradient $\nabla f_L$ are simultaneously small.
We show that if the parameter $\alpha'$ is small enough, then with high probability
the topology of $Z(\widetilde{f}_L)$ depends only on the signs of the eigenvalues
of the Hessian $\nabla^2 f_L(p)$ and the signs of $\widetilde{f}_L(p)$,
$p\in\Cr(\alpha)$. We will describe how these signs determine the structure of the
zero set $Z(\widetilde{f}_L)$ in small neighbourhoods of the points $p\in\Cr(\alpha)$.
Outside these neighbourhoods the zero lines of $Z(\widetilde{f}_L)$ stay close to the
ones of $Z(f_L)$.

First, we consider the critical points $p\in\Cr(\alpha)$ for which both eigenvalues of
the Hessian $\nabla^2 f_L(p)$ have the same sign, i.e., the points
which are local extrema of $f_L$.
In this case, we show that there exists
a disk $D(p, \delta)$ centered at $p$ of a small radius $\delta$ such that
with high probability $ Z(\widetilde{f}_{L})\cap D(p, \delta)$ either consists of a
simple loop encircling the point $p$ when the sign of $\widetilde{f}_{L}(p)$
is opposite to that of the eigenvalues of $\nabla^2 f_L(p)$, or is empty when these
signs coincide. We call such connected components of $Z(\widetilde{f}_L)$ {\em blinking circles}.

\begin{figure}[h]
\includegraphics[width=0.6\textwidth]{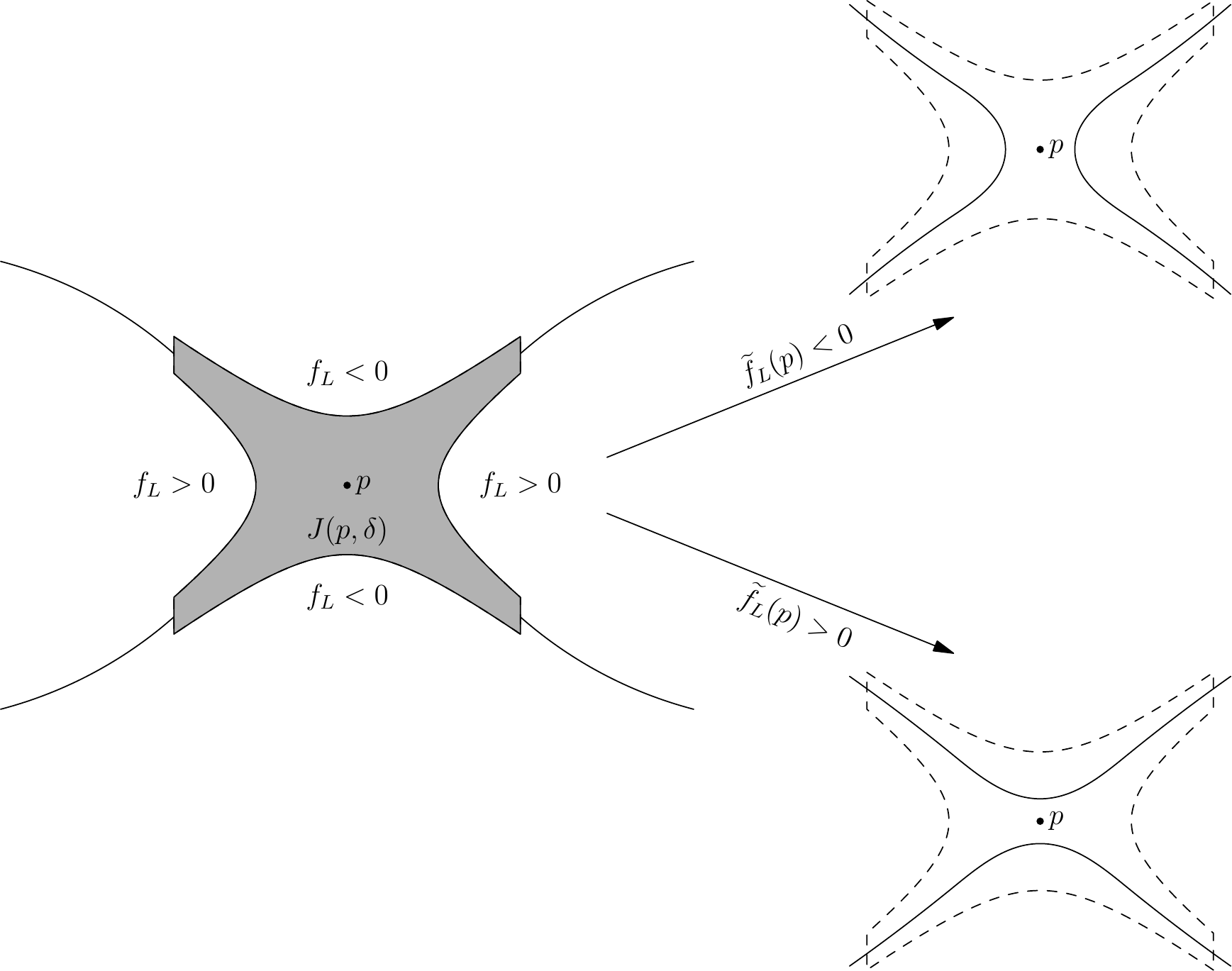}
\caption{The sign of $\widetilde{f}_L$ determines the structure of the zero set $Z(\widetilde{f}_L)$
near the saddle point $p$}
\label{zfigureA}
\end{figure}

Now, we turn to the case when the eigenvalues of the Hessian
$\nabla^2 f_L(p)$ have opposite signs, i.e., to the saddle points of $f_L$.
In this case, the situation is more intricate. We define a degree four
graph $ G(f_L) $ embedded in $\bS^2(L)$. Its vertices are small neighbourhoods
$J(p, \delta)$ of saddle points $p\in\Cr(\alpha)$.
The edges are arcs in the set $ Z(f_L) $ that connect these neighbourhoods.
We will show that with high probability the collection of signs
$\{\operatorname{sgn}(\widetilde{f}_L(p))\colon p\in\Cr(\alpha)\}$ determines how the graph $G(f_L)$ is turned into a collection of loops in $Z(\widetilde{f}_L)$ which we will call {\em the Bogomolny-Schmit loops}. Figure~\ref{zfigureA}
illustrates how the sign of $\widetilde{f}_L$ at the saddle point $p\in\Cr(\alpha)$ determines
the structure of the zero set $Z(\widetilde{f}_L)$ in a small neighbourhood of the saddle point $p\in \Cr(\alpha)$.

We see that in both cases {\em the fluctuations in the number of connected components of
$Z(\widetilde{f}_L)$ are caused by fluctuations in the signs}
\[
\sgn(\widetilde{f}_L(p)) = \sgn(\sqrt{1-\alpha'^2}f_L(p) +\alpha' g_L(p) ),
\quad p\in\Cr (\alpha)\,.
\]

We show that, since the points of the set $\Cr(\alpha)$ are well separated
and the covariance kernel $k_L(d)$ decays at least as a power of $d$, with probability very close to one, we can replace the values $g_L(p)$, $p\in\Cr(\alpha)$, by a collection of
{\em independent} standard Gaussian random variables.

\medskip
Thus, conditioning on $f_L$, we may assume that the values of $\widetilde{f}_L$ at $\Cr(\alpha)$ are
independent normal random variables (not necessarily mean zero).
To conclude, we apply Lemma~\ref{lemma:variance_clusters} on the variance of the number of loops generated by percolation-like processes on planar graphs of degree $4$.

\medskip Note that this chain of arguments can be viewed as the first, although very modest,
step towards justification of the Bogomolny-Schmit heuristics.

\section{Notation}

Throughout the paper, we will be using the following notation:
\begin{itemize}
\item $L$ is a large parameter which tends to $+\infty$. We always assume that $L\ge 1$.
\item $\bS^2(L)$ denotes the sphere in $\bR^3$ centered at
the origin and of radius $L$, while, as usual, $\bS^2$ denotes the unit sphere in $\bR^3$.
By $d_L$ we denote the spherical distance on $\bS^2(L)$.
By $D(x, \rho)\subset \bS^2(L)$ we denote
the open spherical disk of radius $\rho$ centered at $x$.
\item $\overline{G}$ is the closure of the set $G$.
\item We use the abbreviations {\em a.s.} for ``almost surely'',
{\em w.o.p.} for ``with overwhelming probability'', which means
that the property in question holds outside an event of
probability $O(L^{-C})$ with {\em every} $C>0$,
and {\em w.h.p.} for ``with high probability'', which means the property in question holds
outside an event of probability $O(L^{-c})$ with {\em some} $c>0$.
\item $C$ and $c$ (with or without indices)
are positive constants that might only depend on the parameters
in the definition of the regular Gaussian ensemble $f_L$
($\mathcal C^{3+\nu}$-smoothness and the power decay of correlations).
One can think that the constant $C$ is large (in particular, $C\ge 1$),
while the constant $c$ is small (in particular, $c\le 1$).
The values of these constants are irrelevant for our purposes and may vary from line to line.
\item
$A\lesssim B$ means $A \le C\cdot B$, $A \gtrsim B$ means $A \ge c\cdot B$,
and $A\simeq B$ means that $A \lesssim B$ and $A\gtrsim B$ simultaneously.
The sign $\ll$ means ``sufficiently smaller than'', and $\gg$ means
``sufficiently larger than''.
\end{itemize}

\section{Preliminaries}

For the reader's convenience, we collect here standard facts that we will
be using throughout this work.

\subsection{Local coordinates}
It will be  convenient to associate with each point $p\in\bS^2(L)$ its
own coordinate chart. For $p\in\bS^2(L)$, let $\Pi_p$ be a plane in $\bR^3$ passing
through the origin and orthogonal to $p$. The Euclidean structure on $\Pi_p$ is
inherited from $\bR^3$. By $\bS^2_p(L)$ we denote the hemisphere of $\bS^2(L)$ centered
at $p$. By
\[
\Psi_p\colon \bigl\{X\in\Pi_p\colon |X|\le L \bigr\}
\to \bS^2_p(L)
\]
we denote the map inverse to the orthogonal projection. Note that
\[
|X-Y| \le d_L(\Psi_p(X), \Psi_p(Y)) \le 2 |X-Y|
\]
whenever $|X|, |Y|\le \tfrac12 L$.

Let $f\colon \bS^2(L)\to \bR$ be a smooth function. We put $F_p=f\circ \Psi_p$ and
identify ${\rm d}^kf (p)$ with ${\rm d}^k F_p(0)$, i.e., with a $k$-linear form on $\Pi_p$.
Then the gradient $\nabla f(p)= \nabla F_p(0)$ is a vector
in $\Pi_p$ such that $ {\rm d}f(p)(v) = \langle \nabla f(p), v \rangle $, $v\in\Pi_p$,
and the Hessian $H_f(p)=\nabla^2 f(p)=\nabla^2 F_p(0)$ is a self-adjoint operator on $\Pi_p$
such that ${\rm d}^2 f(p)(u, v)=\langle H_f(p)u, v \rangle$, $u, v\in \Pi_p$.

The notation $f\in C^k(\bS^2(L))$ means that,
for every $p\in\bS^2(L)$,
$ f\circ \Psi_p \in C^k\bigl(\{X\in\Pi_p\colon |X|\le \tfrac12 L\} \bigr) $,
and
\[
\| f \|_{C^k} \stackrel{\rm def}= \sum_{j=0}^k \max_{\bS^2(L)} \| {\rm d}^j f \|.
\]
Obviously,
$ \| f \|_{C^k} \le
\max_{p\in \bS^2(L)}\,\| f \circ \Psi_p \|_{C^k(\{|X|\le \frac12 L)\}} $.
In the other direction, it is not difficult to see that
if $f\in C^k(\bS^2(L))$, then
$ \| f \circ \Psi_p \|_{C^k(\{|X|\le \frac12 L\})} \le C_k \| f \|_{C^k} $.

\subsection{Statistical properties of the gradient and the Hessian}

Fix $p\in\bS^2(L)$ and the orthogonal coordinate system $(X_1, X_2)$ on the plane
$\Pi_p$, and set $\partial^k_{i_1\, \ldots \, i_k} f(p)
= \partial^k_{X_{i_1}\, \ldots \, X_{i_k}} F_p(0)$, where, as above, $F_p=f \circ\Psi_p$.

\subsubsection{Independence}

To simplify the notation, next we will deal with the case $L=1$. The general case can be easily obtained by scaling.
\begin{lemma}\label{lemma:independence}
Let $f$ be a $C^{2+\nu}$-smooth random Gaussian function on the
sphere $\bS^2$ whose distribution is invariant with respect to the isometries of the sphere.
Then the following Gaussian random variables are independent:

\smallskip\noindent{\rm (i)}
$f(p)$ and $\nabla f(p)$, as well as $\nabla f(p)$ and $\nabla ^2 f(p)$;

\smallskip\noindent{\rm (ii)}
$\partial_1 f(p)$ and $\partial_2 f(p)$;

\smallskip\noindent{\rm (iii)}
$f(p)$ and $\partial^2_{1, 2} f(p)$;

\smallskip\noindent{\rm (iv)}
$\partial^2_{1, 1} f(p)$ and $\partial^2_{1, 2} f(p)$, as well as
$\partial^2_{2, 2} f(p)$ and $\partial^2_{1, 2} f(p)$.

\end{lemma}

\noindent{\em Proof}: WLOG, we assume that $p$ is the North Pole of the sphere
$\bS^2$ and suppress the dependence on $p$, letting $\Psi=\Psi_p$.
Then,
$$ \Psi (X_1, X_2) = \bigl( X_1, X_2, \sqrt{1-(X_1^2+X_2^2)}\, \bigr), $$
and $F=f\circ \Psi$ is a Gaussian function on the unit disk
$\{|X_1|^2+|X_2|^2<1\}$ with the covariance
\begin{align*}
K(X, Y) &= \bE [F(X) F(Y)] \\
&= k\bigl( d(\Psi (X), \Psi (Y))\, \bigr) \qquad (d \text{\ is\ the\ spherical\ distance}) \\
&= k\bigl( \arccos\langle \Psi (X), \Psi (Y)\rangle \bigr)
\\
&= k\bigl( \arccos \bigl( \sum_{j=1}^2 X_j Y_j +
\bigl( 1-\sum_{j=1}^2 X_j^2\bigr)^{\frac12} \bigl( 1-\sum_{j=1}^2 Y_j^2\bigr)^{\frac12}\bigr) \bigr),
\end{align*}
Note that
\begin{equation}\label{eq:K}
K(-X_1, X_2, -Y_1, Y_2) = K(X_1, -X_2, Y_1, -Y_2) = K(X_1, X_2, Y_1, Y_2).
\end{equation}

To prove properties (i)--(iv), we need to check the corresponding statements for $F(0)$, $\partial_{X_j} F(0)$,
and $\partial^2_{X_i X_j} F(0)$. The covariances of these random variables can be computed using the relations
\[
\bE \Bigl[ \frac{\partial^\ell F(X)}{\partial^{\ell_1}X_1 \partial^{\ell_2} X_2} \,
\frac{\partial^m F(Y)}{\partial^{m_1}Y_1 \partial^{m_2} Y_2}
\Bigr] \, =  \,
\frac{\partial^{\ell+m} K(X, Y)}
{\partial^{\ell_1} X_1 \partial^{\ell_2} X_2 \partial^{m_1}Y_1 \partial^{m_2} Y_2}\,.
\]
The rest follows by differentiation of relations~\eqref{eq:K}.
For instance,
\[
\frac{\partial^2 K}{\partial Y_1 \partial Y_2} (X_1, X_2, Y_1, Y_2)
= - \frac{\partial^2 K}{\partial Y_1 \partial Y_2} (-X_1, X_2, -Y_1, Y_2),
\]
whence,
\[
\bE\bigl[ F(0) \partial^2_{Y_1 Y_2} F(0) \bigr] = \frac{\partial^2 K}{\partial Y_1
\partial Y_2} (0, 0) = 0\,.
\]
Similarly,
\[
\frac{\partial^4 K}{\partial^2 X_1 \partial Y_1 \partial Y_2} (X_1, X_2, Y_1, Y_2)
= - \frac{\partial^4 K}{\partial^2 X_1 \partial Y_1 \partial Y_2} (-X_1, X_2, -Y_1, Y_2),
\]
whence,
\[
\bE\bigl[ \partial^2_{X_1}F(0) \partial^2_{Y_1 Y_2} F(0) \bigr]
= \frac{\partial^4 K}{\partial^2_{X_1} \partial Y_1 \partial Y_2} (0, 0) = 0\,,
\]
and so on \ldots\,  . \mbox{} \hfill $\Box$

\subsubsection{Non-degeneracy}

\begin{lemma}\label{lemma:non-degeneracy}
Let $(f_L)$ be a regular Gaussian ensemble, and
let $p\in\bS^2(L)$. Then
\[
\liminf_{L\to\infty}\, \bE\bigl[ (\partial_{j} f_L(p))^2 \bigr] >0, \qquad j\in\{1, 2\},
\]
and
\[
\liminf_{L\to\infty}\, \bE\bigl[ (\partial^2_{i, j} f_L(p))^2 \bigr] >0, \qquad i,j\in\{1, 2\}.
\]
\end{lemma}

\noindent{\em Proof}: Again, we assume that $p$ is
the North Pole of the sphere $\bS^2(L)$. Put
$ \Psi_L (X_1, X_2) = \bigl( X_1, X_2, \sqrt{L^2-(X_1^2+X_2^2)}\, \bigr)$,
and consider the Gaussian functions
$F_L = f_L \circ \Psi_L$ defined in the disks $\tfrac12 L \bD$. The corresponding
covariances
$ \bE[F_L(X)F_L(Y)] = K_L(X, Y) $
are $C^{3+\nu, 3+\nu}$-smooth on
$\tfrac12 L\bD \times \tfrac12 L \bD $ with some $\nu>0$. Their partial derivatives
up to the third order are bounded locally uniformly in $L\ge L_0$.
Hence, by a version of the Arzel\'a-Ascoli theorem, any sequence $K_{L_j}$
contains a locally uniformly $C^{2+\nu, 2+\nu}$-convergent subsequence.

The limiting function $K$ is a $C^{2+\nu, 2+\nu}$-smooth Hermitean-positive
function on $\bR^2~\times~\bR^2$ which depends only on the Euclidean distance $|X-Y|$.
Hence, by Bochner's theorem, it is a Fourier integral of a finite positive
rotation-invariant measure $\rho$:
\[
K(X, Y) = \int_{\bR^2} e^{2\pi{\rm i}\langle \la, X-Y \rangle}\, {\rm d}\rho(\la),
\]
where
\[
\int_{\bR^2} |\la|^4\, {\rm d}\rho(\la) <\infty\,.
\]

Furthermore, since
\[
\bigl| K_L(X, Y) \bigr| \lesssim \bigl( 1+d_L(\Psi_L(X), \Psi_L(Y) )\bigr)^{-\gamma},
\]
uniformly in $L\ge L_0$, the limiting function $K$ satisfies
\[
\bigl| K(X, Y) \bigr| \lesssim \bigl( 1 + |X-Y| \bigr)^{-\gamma}.
\]
Therefore, the measure $\rho$ cannot degenerate to the point measure
at the origin.

The rest is straightforward.
Suppose, for instance, that for some sequence $L_j\to\infty$,
\[
\lim_{j\to\infty} \bE\bigl[ (\partial_{1} f_{L_j} (p))^2\bigr] = 0.
\]
Then
\[
\lim_{j\to\infty} \frac{\partial^2 K_{L_j}}{\partial X_1 \partial Y_1} (0, 0) = 0. \]
Passing to a subsequence, we conclude that
\[
(2\pi {\rm i})^2\, \int_{\bR^2} \la_1^2\, {\rm d}\rho(\la)
= -\frac{\partial^2}{\partial X_1 \partial Y_1}\,
\int_{\bR^2} e^{2\pi{\rm i}\langle \la, X-Y \rangle}\, {\rm d}\rho(\la)\, \Big|_{X=Y=0}
= 0\,.
\]
Since the measure $\rho$ is positive and rotation-invariant, this is possible only when
$\rho$ is a point mass at the origin. This contradiction concludes the proof.
\mbox{} \hfill $\Box$

\subsubsection{Power decay of correlations}

The power decay of the correlations between $f_L(p)$ and $f_L(q)$ when $d_L(p, q)$
is large and the a priori $C^{3,3}$-smoothness of the covariance yield the power decay of
correlations between the Gaussian vectors
\[
v(p) = (f_L(p), \nabla f_L(p)) = (F_p(0), \nabla F_p(0)) =
(F_p(0), \partial_{X_1} F_p(0), \partial_{X_2} F_p(0) )
\]
and
\[
v(q) = (f_L(q), \nabla f_L(q)) = (F_q(0), \nabla F_q(0)) =
(F_q(0), \partial_{Y_1} F_q(0), \partial_{Y_2} F_q(0) )
\]
(we keep fixed the coordinate systems $(X_1, X_2)$ and $Y_1, Y_2)$
in the planes $\Pi_p$ and $\Pi_q$).

\begin{lemma}\label{lemma:correl-decay}
Let $(f_L)$ be a regular Gaussian ensemble.
Then, for any $p, q \in\bS^2(L)$,
\[
\max_{1\le i, j \le 3}\ \bigl| \bE[v_i(p) v_j(q)]\bigr|
\lesssim \bigl(1+d_L(p, q) \bigr)^{- \gamma/4}.
\]
\end{lemma}

\noindent{\em Proof}:
Put $K_{p, q}(X, Y) = \bE\bigl[F_p(X) F_q(Y)\bigr]$, where
$F_p = f_L \circ \Psi_p$, and $(X, Y)\in \overline{Q} \times \overline{Q}$,
where $\overline{Q}=[-1, 1]\times [-1, 1]$.
We have $\bE [F_p(0) \partial_{Y_i} F_q(0)]  = \partial_{Y_i} K_{p, q}(0, 0) $, and
$ \bE [\partial_{X_i} F_p(0) \partial_{Y_j} F_q(0)] =
\partial^2_{X_i Y_j} K_{p, q}(0, 0) $.

There is nothing to prove if $d_L(p, q) \le 1$, so we assume that $d_L(p, q) \ge 1$.
Then, $\| K_{p, q}\|_{C(\overline{Q}\times\overline{Q})}\!\lesssim\!d^{-\gamma}$ and
$\| K_{p, q}\|_{C^{2, 2}(\overline{Q}\times\overline{Q})} \lesssim 1$.
Then, by the classical Landau-Hadamard inequality\footnote{We use it in the following
form. If $h\colon [0, 1] \to \bR$ is a $C^2$-smooth function and
$M_j = \max_{[0, 1]} |h^{(j)}|$, $0 \le j \le 2$, then $M_1 \lesssim \max(M_0, \sqrt{M_0 M_2}\,)$}
applied to the functions $X_i \mapsto K_{p, q}(X, Y)$ and
$Y_j \mapsto K_{p, q}(X, Y)$, we have
$ \| \partial_{X_i} K_{p, q}\|_{C(\overline{Q}\times\overline{Q})}
\lesssim d^{-\gamma/2}$ and
$\| \partial_{Y_j} K_{p, q}\|_{C(\overline{Q}\times\overline{Q})} \lesssim d^{-\gamma/2}$,
$i, j = 1, 2$.
Applying the Landau-Hadamard inequality again, this time to the functions
$Y_j \mapsto \partial_{X_i} K_{p, q}(X, Y)$, we get
$ \| \partial^2_{X_i Y_j} K_{p, q}\|_{C(\overline{Q}\times\overline{Q})}
\lesssim d^{-\gamma/4}$, $i, j = 1, 2$.

In particular, these estimates hold at $X=Y=0$,
which gives us what we needed.  \mbox{} \hfill $\Box$

\begin{itemize}
\item
{\em To simplify our notation, in what follows, we assume that the parameter $\gamma>0$ is chosen so that the correlations between the Gaussian vectors $(f_L(p), \nabla f_L(p))$ and $(f_L(q), \nabla f_L(q))$ decay as $(1+d_L(p, q))^{-\gamma}$}.
\end{itemize}

\subsection{A priori smoothness of $f_L$}

Let $(f_L)$ be a regular Gaussian ensemble.
Quite often, we will be using the following a priori bound:
\begin{itemize}
\item {\em w.o.p, $ \| f_L \|_{C^3(\bS^2(L))} < \log L $}.
\end{itemize}
This bound immediately follows from the classical estimate
\[ \bP\bigl\{ \| f_L \|_{C^3(\bS^2(L))} > t \bigr\} \le CL^2 e^{-ct^2}. \]
For a self-contained proof see, for instance,~\cite[Sections~A9--A11]{NS2}.

\section{Smooth functions with controlled topology of the zero set}

Here we introduce the (non-random) class $C^3(A, \Delta, \alpha, \beta)$
of smooth functions $f$ on $\bS^2(L)$ such that
the number of connected components of the zero set of a small perturbation $\widetilde{f}$ of $f$ can be recovered from the values of $\widetilde{f}$ at the critical points of $f$ with small critical values (provided that the values of $\widetilde{f}$ at these points are not too small). Later, we will show that w.h.p. our random function $f_L$ belongs to this class.

\subsection{The sets $\Cr(\alpha)$, $\Cr(\alpha, \beta)$, and
$\Cr(\alpha, \beta, \Delta)$, and the class $C^3(A, \Delta, \alpha, \beta)$}
\label{subsect:sets_Cr_class_C3}

Given $\alpha, \beta \le 1$ and $\Delta \ge 1$, we let
\[
\Cr(\alpha) = \{p\in \bS^2(L)\colon |f(p)|\le \alpha,\ \nabla f(p) = 0 \},
\]
\[
\Cr(\alpha, \beta) = \{p\in \bS^2(L)\colon |f(p)|\le \alpha, \ |\nabla f(p)| \le\beta\},
\]
and
\[
\Cr(\alpha, \beta, \Delta) = \{p\in \bS^2(L)\colon
|f(p)|\le \alpha, \ |\nabla f(p)| \le\beta, \ \|(\nabla^2 f(p))^{-1}\|_{\tt op} \le \Delta \},
\]
where $\|\,.\,\|_{\tt op}$ stands for the operator norm.

By $C^3(A)$ we denote the class of $C^3$-smooth functions $f$ on $\bS^2(L)$
with $\|f\|_{C^3} \le A$.
Given the parameters
\[\alpha \ll \beta \ll 1 \ll A \ll \Delta, \]
by $C^3(A, \Delta, \alpha, \beta)$ we denote
the class of functions $f\in C^3(A)$ for which
$\Cr(\alpha, \beta) = \Cr(\alpha, \beta, \Delta)$, i.e.,
the Hessian of $f$ does not degenerate ($\|(\nabla^2 f)^{-1}\|_{\tt op} \le \Delta$)
on the almost singular set $\Cr(\alpha, \beta)$ where $f$ and $\nabla f$ are simultaneously
small.

Given $f\in C^3(A, \Delta, \alpha, \beta)$, $\alpha'\ll\alpha$, and
$g\in C^3(A)$, we set
\[
f_t= f+tg, \quad 0\le t \le \alpha'.
\]

Next, we develop a little caricature of the quantitative
Morse theory, which shows
that the collection of signs of $f_t$ at $\Cr(\alpha)$ defines the topology
of the zero set $Z(f_t)$, provided that $\min_{\Cr (\alpha)}|f_t|$ is not too small, and
gives ``an explicit formula'' that recovers the number of connected components of $Z(f_t)$
from this collection of signs and the structure of $Z(f)$.

\subsection{Near any almost singular point there is a unique critical point of $f$}

\begin{lemma}\label{lemma:crit-almost_sing}
Suppose that $f\in C^3(A)$ and $p\in\Cr(\alpha, \beta, \Delta)$
with
\[
1 \ll A \ll \Delta, \quad A \Delta^2 \beta \ll 1.
\]
Then,

\smallskip\noindent{\rm (A)}
the spherical disk $D(p, 2\Delta\beta)$ contains a unique critical
point $z$ of $f$;

\smallskip\noindent{\rm (B)}
there are no other critical points of $f$ in the disk
$D\bigl( p, c (A\Delta)^{-1} \bigr)$;

\smallskip\noindent{\rm (C)} $|f(z)| \le 2\alpha $,
provided that $A\Delta^2\beta^2\ll \alpha$.
\end{lemma}

\medskip\noindent{\em Proof}: We will work on the plane $\Pi_p$, and let $F=f\circ \Psi_p$,
and $\cH_F =\nabla^2 F$. To find the critical point $z$,
we use a simplified Newton's method:
\[
X_{n+1} = X_n - \cH_F(0)^{-1} \nabla F(X_n), \quad X_0=0\,.
\]
Put $\Phi(X)=X-\cH_F(0)^{-1} \nabla F (X)$. First, we check that in the disk
$D\bigl( 0, c (A\Delta)^{-1} \bigr)$ the map $\Phi$ is a $\tfrac13$-contraction.

Indeed,
\begin{multline*}
\Phi(X)-\Phi(Y) = (X-Y) - \cH_F(0)^{-1}\bigl( \nabla F(X)-\nabla F(Y) \bigr) \\
= (X-Y) - \cH_F(0)^{-1}\bigl[ \cH_F(X)(X-Y) + O(A|X-Y|^2) \bigr],
\end{multline*}
whence,
\[
| \Phi(X)-\Phi(Y) | \le \| I-\cH_F(0)^{-1}\cH_F(X) \|_{\tt op}\,
|X-Y| +  \Delta \cdot O(A|X-Y|^2).
\]
Furthermore,
\[
\| I-\cH_F(0)^{-1}\cH_F(X) \|_{\tt op}
\le \| \cH_F(0)^{-1} \|_{\tt op} \, \| \cH_F(0) - \cH_F(X) \|_{\tt op}
\le \Delta \cdot O(A|X|),
\]
and finally,
\[
| \Phi(X)-\Phi(Y) | \le A\Delta \bigl( O(|X|) + O(|X-Y|) \bigr)\, |X-Y|
\le \tfrac13\, |X-Y|,
\]
provided that $|X|, |Y|\le c(A\Delta)^{-1}$ with sufficiently small positive $c$.

Next, we check that the map $\Phi$ preserves the disk
$D(0, c (A\Delta)^{-1})$. We have
\begin{multline*}
\Phi(X) = X- \cH_F(0)^{-1} \nabla F(X) \\
= -\cH_F(0)^{-1} \bigl( \nabla F(X) - \cH_F(0) X \bigr)
= -\cH_F(0)^{-1} \bigl( \nabla F(0) + O(A|X|^2) \bigr),
\end{multline*}
and then,
\[
|\Phi (X)| \le \Delta\, \Bigl( \beta + A \bigl( \frac{c}{A\Delta} \bigr)^2\, O(1) \Bigr)
= \beta\Delta + \frac{c^2}{A\Delta}\, O(1)
< \frac{c}{A\Delta},
\]
provided that $\beta A \Delta^2 \ll 1$ and that the positive constant $c$ is
sufficiently small.

Thus, $\Phi$ has a unique fixed point $Z$ in the disk $D(0, c (A\Delta)^{-1} )$,
and
\begin{multline*}
|Z| = |0-Z|
\le \sum_{n\ge 0} |X_{n}-X_{n+1}|
\le \sum_{n\ge 0} 3^{-n} |X_0-X_1| \\
= \tfrac32\, |0-\Phi(0)|
= \tfrac32\, |\Phi(0)|
= \tfrac32\, |\cH_F(0)^{-1} \nabla F(0)|
\le \tfrac32\, \Delta \beta.
\end{multline*}
At last,
\[
|F(Z)| \le |F(0)| + |\nabla F(0)|\, |Z| + O(A|Z|^2)
\le \alpha + \beta \cdot 2\beta\Delta + O(A \beta^2\Delta^2)
\le 2 \alpha,
\]
provided that $ A \beta^2\Delta^2 \ll \alpha $. \mbox{} \hfill $\Box$

\subsection{Near any point $p\in\Cr(\alpha)$ there is a unique critical point
of $f_t$}

\begin{lemma}\label{LemmaB}
Let $f, g\in C^3(A)$, $f_t=f+tg$ with $0\le t \le \alpha'$.
Suppose that $p\in \Cr(\alpha)$ with $\| \cH_f(p)^{-1}\|_{\tt op} \le \Delta$,
and that
\[
1 \ll A \ll \Delta, \quad
A\alpha' \ll \alpha, \quad A\Delta^2 \alpha \ll 1.
\]
Then there exists a unique critical
point $p_t$ of $f_t$ such that $ d_L(p, p_t) \ll \Delta\alpha $
and $|f_t(p)-f_t(p_t)| \ll A(\Delta \alpha)^2$.
Moreover, there are no other critical points of $f_t$ at distance $\le c (A\Delta)^{-1}$
from $p$.
\end{lemma}

\noindent{\em Proof of Lemma~\ref{LemmaB}}: First, we note that $f_t\in C^3(2A)$
and that $|f_t(p)|\le |f(p)|+A\alpha'$ and $|\nabla f_t(p)|\lesssim A\alpha' \ll \alpha$.
Furthermore, $\| \cH_{f_t}(p)^{-1} \|_{\tt op} \le 2\Delta$. Indeed,
we have
\begin{multline*}
\| \cH_{f_t} (p)^{-1} \|_{\tt op} =
\| \cH_f(p)^{-1} (I + (\cH_{f_t}(p) - \cH_f(p))\cH_f(p)^{-1} )^{-1}\|_{\tt op} \\
\le
\Delta\, \| (I + (\cH_{f_t}(p) - \cH_f(p))\cH_f(p)^{-1} )^{-1} \|_{\tt op}\,.
\end{multline*}
Noting that $ \| \cH_{f_t}(p) - \cH_f(p) \|_{\tt op} \lesssim A\alpha' $, we see that
\[ \| (\cH_{f_t}(p) - \cH_f(p)) \cH_f(p)^{-1} \|_{\tt op} \lesssim \Delta \cdot A\alpha'
\ll \Delta\alpha \ll 1.\]
Therefore, \[ \| (I + (\cH_{f_t}(p) - \cH_f(p))\cH_f(p)^{-1} )^{-1} \|_{\tt op} < 2,\]
and, finally, $\| \cH_{f_t} (p)^{-1} \|_{\tt op} \le 2\Delta$.

Then, by Lemma~\ref{lemma:crit-almost_sing} (applied to the function $f_t$ with $\beta=A\alpha'$),
there exists a unique critical point $p_t$ of $f_t$ with
\[
d_L(p, p_t) \le 2 (A\alpha') \cdot 2\Delta \ll \Delta\alpha.
\]
Besides, for $d_L(p, x)\ll\Delta\alpha$, we have
$ d_L(p_t, x) \ll \Delta \alpha $ and then
$|\nabla f_t(x)| = |\nabla f_t(x) - \nabla f_t(p_t) | \ll A\cdot \Delta\alpha$,
whence, \[ |f_t(p)-f_t(p_t)| \ll A(\Delta\alpha)^2. \]

At last, by part~B of Lemma~\ref{lemma:crit-almost_sing}, there are no other critical points of $f_t$ at distance $\le c(A\Delta)^{-1}$ from $p$. \mbox{} \hfill $\Box$

\subsection{Local matters}

Given $f\in C^3(A)$, $p\in\Cr(\alpha)$, $\| \cH_f(p)^{-1} \|_{\tt op} \le\Delta$,
we look at the behaviour of $f_t$ in the $\delta$-neighbourhood of $p$. As above, $f_t=f+tg$,
with $g\in C^3(A)$, and $0\le t \le \alpha'$. Throughout this section we assume that the
parameters $\alpha'$, $\alpha$, $\delta$, $A$ and $\Delta$ satisfy the following set
of conditions
\begin{equation}\label{eq:parameters}
\alpha \ll 1 \ll A \ll \Delta,
\quad A\alpha' \ll \alpha \ll A^{-2} \Delta^{-3},
\end{equation}
and that
\begin{equation}\label{eq:delta}
\delta = c(A \Delta)^{-1}
\end{equation}
with sufficiently small constant $c$. Note that these conditions are more restrictive
than the ones used in Lemma~\ref{LemmaB}, so we will be using freely that lemma.

\subsubsection{Local extrema}\label{subsect:blinking-circles}
First, we consider the case when the Hessian $\cH_f(p)$ is positive or negative
definite, that is, its eigenvalues have the same sign.
With a little abuse of terminology, we say that the function $f_t$ is convex (concave) in
$D(p, \delta)$ if the function $f_t\circ \Psi_p$ is convex (correspondingly, concave)
in $\Psi_p^{-1} D(p, \delta)\subset \Pi_p$.

\begin{lemma}\label{LemmaC}
Suppose that the eigenvalues of the Hessian $\cH_f(p)$ have the same sign
and that conditions~\eqref{eq:parameters} and~\eqref{eq:delta} hold. Then

\smallskip\noindent{\rm (i)}
the function $f_t$ is either concave or convex function in $D(p, \delta)$,

\smallskip\noindent{\rm (ii)}
the function $f_t$ does not vanish on $\partial D(p, \delta)$, and moreover, the sign
of $f_t\big|_{\partial D(p, \delta)}$ coincides with the sign of the eigenvalues
of $\cH_f(p)$.

\end{lemma}

\noindent{\em Proof of Lemma~\ref{LemmaC}}:
Put $F=f\circ \Psi_p$, $F_t = f_t\circ \Psi_p$,
and suppose, for instance, that $\cH_F(0)=\cH_f(p)$ is positive
definite (otherwise, replace $f$ by $-f$), that is,
$\langle \cH_F(0) x, x \rangle \ge \Delta^{-1} |x|^2 $.
Then, for any $X\in \Psi_p^{-1}D(p, \delta)$, we have
\begin{multline*}
\langle \cH_F(X) x, x \rangle
\ge \langle \cH_F(0) x, x \rangle - \| \cH_F(X)-\cH_F(0)\|_{\tt op} |x|^2 \\
\ge (\Delta^{-1} - CA|X|) |x|^2 \ge (2\Delta)^{-1} |x|^2
\end{multline*}
since $\| \cH_F(X)-\cH_F(0)\|_{\tt op}\lesssim A|X| $ and
$|X| \le \delta = c(A\Delta)^{-1}$
with sufficiently small $c$. Noting that
$ \| \cH_{F_t}(X)-\cH_{F}(X) \|_{\tt op} \lesssim A\alpha' \ll \Delta^{-1} $, we
get $ \langle \cH_{F_t}(X) x, x \rangle \ge (4\Delta)^{-1} |x|^2 $ which proves (i).

To prove (ii), we take $X\in \Psi_p^{-1} \partial D(p, \delta)$.
Then
\begin{multline*}
F_t(X) = F_t(0) + \langle \nabla F_t(0), X \rangle + \tfrac12\,
\langle \cH_{F_t}(0) X, X \rangle
+ O(A\delta^3) \\
\ge \tfrac12\, \langle \cH_{F_t}(0) X, X \rangle
-  |F_t(0)| - | \langle \nabla F_t(0), X \rangle | - O(A\delta^3)\,.
\end{multline*}
Furthermore, using that $|F_t(0)| = |f_t(p)| \le \alpha + O(A\alpha') \lesssim \alpha$
and that $ | \langle \nabla F_t(0), X \rangle | \le |\nabla F_t(0)| \cdot |X|
= O(\alpha' A) \cdot \delta \ll \alpha$, we conclude that
\begin{multline*}
F_t(X)\, \stackrel{|X|\ge \delta/2}\ge\, \tfrac12\,
(4\Delta)^{-1} (\delta/2)^2 - O(\alpha + A\delta^3)
= (32\Delta)^{-1} \delta^2 - O(A\delta^3) \\
= \bigl( (32\Delta)^{-1}  - O(A\delta) \bigr)\, \delta^2
\, \stackrel{\delta= c(A\Delta)^{-1}}\ge\,
\bigl( \tfrac1{32} - C\cdot c  \bigr)\, \Delta^{-1} \delta^2
\ge \tfrac1{64} \Delta^{-1} \delta^2\,,
\end{multline*}
provided that the constant $c$ in~\eqref{eq:delta} (the definition of $\delta$)
was chosen so small that $C\cdot c \le \tfrac1{64}$. \mbox{} \hfill $\Box$

\medskip\noindent\underline{Summary}: Let $p$ be a local extremum of $f$. Suppose that
conditions~\eqref{eq:parameters} and~\eqref{eq:delta} hold.
\begin{itemize}
\item Then, $Z(f_t)\cap \partial D(p, \delta)= \emptyset$,
\item $Z(f_t)\cap D(p, \delta)$ is either empty, or homeomorphic to $\bS^1$, or
a singleton.
\item Suppose that $|f_t(p)| \gtrsim A (\Delta\alpha)^2$. Then, by Lemma~\ref{LemmaB},
$f_t(p_t)$ has the same sign as $f_t(p)$. Therefore,
$Z(f_t)\cap D(p, \delta)=\emptyset$ whenever $f_t(p)$ and the eigenvalues of $\cH_f(p)$
have the same sign, and $Z(f_t)\cap D(p, \delta)$ is homeomorphic to $\bS^1$
whenever $f_t(p)$ and the eigenvalues of $\cH_f(p)$ have opposite signs.
\end{itemize}

\subsubsection{Saddle points}\label{subsect:saddle-points}
Now, we turn to the case when $p\in\Cr(\alpha)$ is a saddle point of $f$, that is, the eigenvalues of $\cH_f(p)$ have opposite signs.
We will work on the plane $\Pi_p$ and set $F=f\circ \Psi_p$, $G=g\circ \Psi_p$,
$F_t = f_t\circ\Psi_p = F + tG$. By
$ H(X)=\langle \cH_F(0)X, X \rangle $
we denote the quadratic form generated by the Hessian $\cH_F(0)$. WLOG, we assume
that
\[
H(X)=aX_1^2-bX_2^2, \quad \Delta^{-1} \le a \le b \le A.
\]
We take $\delta=c(A\Delta)^{-1}$ with a sufficiently small positive constant $c$,
set
\[
J(\delta) \stackrel{\rm def}=
\bigl\{|H|\le a\delta^2 \bigr\} \bigcap \bigl\{ |X_1|\le 3\delta \bigr\}
\]
and call this set {\em a joint}.
\begin{figure}[h]
\includegraphics[width=0.6\textwidth]{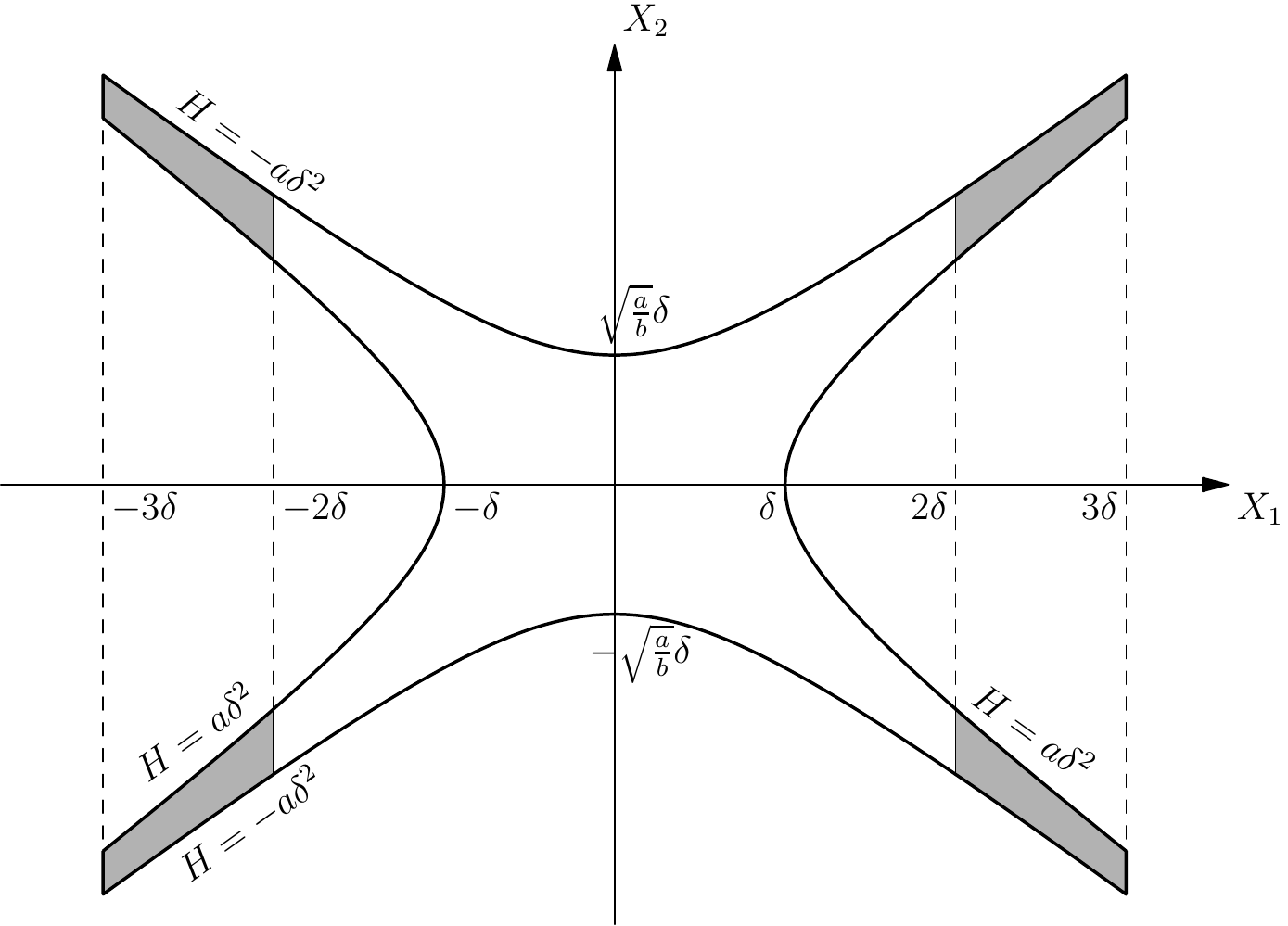}
\caption{Joint $J(\delta)$ with four terminals}
\label{zfigure1}
\end{figure}
By
\[
\partial^* J(\delta) \stackrel{\rm def}= \bigl\{|H|= a\delta^2 \bigr\} \bigcap \bigl\{ |X_1|\le 3\delta \bigr\}
\]
we denote the curvilinear part of the full boundary
$\partial J(\delta)$ of the joint $J(\delta)$.

\begin{lemma}\label{LemmaE}
Suppose that the eigenvalues of $\cH_f(p)$ have the opposite signs and that
conditions~\eqref{eq:parameters} and~\eqref{eq:delta} hold.
Then the function $F_t$ does not vanish on
$ \partial^* J(\delta) $.
Moreover, the signs of $F_t$ and $H$ coincide on $\partial^* J(\delta)$.
\end{lemma}

\noindent{\em Proof of Lemma~\ref{LemmaE}}:
Everywhere in $J(\delta)$ we have
\[
F_t = F(0) + \tfrac12\, H + O(A\delta^3) + \alpha' G
= \tfrac12\, H + O\bigl( \alpha + (\delta^3 + \alpha') A \bigr)
\, \stackrel{\alpha' \ll\delta^3}=\, \tfrac12 H + O\bigl( A\delta^3 \bigr).
\]
Furthermore, on $\partial^* J(\delta)$ we have
$ |H|=a\delta^2 \ge c^2 A^{-2}\Delta^{-3} $,
while $A\delta^3 = c^3 A^{-2}\Delta^{-3}$.
This proves the lemma. \mbox{} \hfill $\Box$

\medskip The set
$\bigl\{|H|\le a\delta^2 \bigr\} \bigcap \bigl\{ 2\delta \le |X_1|\le 3\delta \bigr\}$
consists of $4$ disjoint curvilinear quadrangles. We call them {\em  terminals}
and denote them by $T_i$, $1\le i \le 4$.

\begin{lemma}\label{LemmaF}
Under the same assumptions as in Lemma~\ref{LemmaE},
each of the sets $Z(F_t)\cap T_i$ consists of one curve which joins the vertical segments
on the boundary of $T_i$.
\end{lemma}

\noindent{\em Proof of Lemma~\ref{LemmaF}}:
Everywhere in $J(\delta)$ we have
\[
(F_t)_{X_2} = \tfrac12\, H_{X_2} + (F_{X_2} - \tfrac12\, H_{X_2}) + tG_{X_2}.
\]
Hence,
\[
|(F_t)_{X_2}| \ge b |X_2| - CA|X|^2 - O(A\alpha')\,.
\]
In each of the terminals $T_i$,
\[
b|X_2| \ge b\cdot \sqrt{3\frac{a}{b}}\cdot \delta > \sqrt{ab}\cdot \delta
\ge \frac{\delta}{\Delta}
\]
(in the last estimate we use that $b\ge a\ge \Delta^{-1}$),
and
\[
|X|^2 \le (3\delta)^2 + \Bigl( \sqrt{10\,\frac{a}{b}}\cdot \delta \Bigr)^2 \le 19 \delta^2,
\]
whence
\[
C A|X|^2 \le 19C A\delta^2 = 19C A\,\frac{c}{A\Delta}\cdot \delta \le \frac{\delta}{2\Delta}\,,
\]
provided that the constant $c$ in the definition of $\delta$ is sufficiently small.
Furthermore, $A\alpha'$ is also much smaller than $\Delta^{-1} \delta$
(since $A\alpha' \ll A^{-1}\Delta^{-2}$). Thus, $|(F_t)_{X_2}|>0$ everywhere
in $T_i$. It remains to recall that, by Lemma~\ref{LemmaE}, the function $F_t$
has at least one change of sign on each vertical section of
$T_i$. Therefore, by the implicit function theorem,
$Z(F_t)\cap T_i$ is a graph of a smooth function. \mbox{} \hfill $\Box$

\medskip
Under the same assumptions as in Lemmas~\ref{LemmaE} and~\ref{LemmaF},
by Lemma~\ref{LemmaB}, the joint $J(\delta)$ contains only one critical  point
$X^t=(X^t_1, X^t_2)$ of $F_t$, and $|X^t|\ll \Delta\alpha$.
Consider the sets
\[
I_1 = \bigl\{X=(X_1, X^t_2)\colon X_1\in\bR\bigr\} \cap J(\delta),
\quad I_2 = \bigl\{X=(X^t_1, X_2)\colon X_2\in\bR\bigr\} \cap J(\delta).
\]
Since
\begin{multline*}
|X^t| \ll \Delta\alpha \ll (A\Delta)^{-2} = c^{-1} (A\Delta)^{-1} \delta
\ll (A\Delta)^{-1/2} \delta \\
\le \sqrt{a/b} \cdot \delta \qquad (\Delta^{-1} \le a \le b \le A),
\end{multline*}
it is easy to see that both sets are the segments.

\begin{lemma}\label{LemmaG}
Under the same assumptions as in Lemmas~\ref{LemmaE} and~\ref{LemmaF},
the only extremum of the restriction of the function $F_t$ to the segment
$I_1$ is a local minimum at $X_1=X^t_1$, and the only extremum of the restriction
of the function $F_t$ to the segment $I_2$ is a local maximum at $X_2=X^t_2$.
\end{lemma}

\noindent{\em Proof of Lemma~\ref{LemmaG}}: Consider the function
$X_1 \mapsto F_t(X_1, X_2^t)$. It has a critical point at $X_1=X_1^t$,
the second derivative at this point is not less than
\[
2a - CA\alpha' - CA |X^t|
\ge \frac2{\Delta} - CA(\alpha'+\Delta \alpha) \ge \frac1{\Delta}
\]
(since $A\Delta^2 \alpha \ll 1$, and $A\alpha' \ll \alpha$),
and the  $C^3$-norm of $F_t$ is bounded by $CA$. Therefore, for $X_1>X_1^t$, we have
\[
(F_t)_{X_1}(X_1, X_2^t) \ge 2a (X_1-X^t_1) - CA (X_1-X^t_1)^2
\ge \Delta^{-1} (X_1-X^t_1) - CA|X_1-X^t_1|^2.
\]
Note that RHS of the last expression is positive since
\[
X_1-X^t_1 < 2\delta = 2c (A\Delta)^{-1}
\]
with sufficiently small constant $c$. Similarly, $(F_t)_{X_1}(X_1, X_2^t) < 0$ for
$X_1<X_1^t$.

The proof of the second statement is almost identical
and we skip it. \hfill $\Box$

\begin{lemma}\label{LemmaH}
Suppose that $F_t(X^t)\ne 0$ (i.e., zero is not a critical value of the
restriction of the function $F_t$ to the joint $J(\delta)$).
Then, under the same assumptions as in Lemmas~\ref{LemmaE}, ~\ref{LemmaF}
and~\ref{LemmaG}, the set $Z(F_t)\cap J(\delta)$ consists of two connected components,
which enter and exit the joint $J(\delta)$ through the terminals $T_i$.

Furthermore, the set $\{F_t\ne 0\}\cap J(\delta)$ consists of three connected components.
One of them contains $F_t(X^t)$, while on the other two components $F_t$ has the sign opposite
to the sign of $F_t(X^t)$.
\end{lemma}

\noindent{\em Proof of Lemma~\ref{LemmaH}}:
Since zero is not a critical value of the restriction $F_t\big|_{J(\delta)}$, the set
$Z(F_t)\cap J(\delta)$ consists of a finitely many disjoint smooth curves.
By Lemma~\ref{LemmaF}, this set
has at least two connected components, the ones that
enter and exit the joint $J(\delta)$ through the terminals. If there exists a third
component, then, again by Lemma~\ref{LemmaF}, it cannot intersect the terminals,
while, by Lemma~\ref{LemmaE}, it also cannot intersect the rest of the boundary
$\partial^* J(\delta)$. Hence, it stays inside the joint.
Therefore, it is a closed curve which bounds a domain $G$ with $\bar G \subset J(\delta)$. Since $F_t$ vanishes on $\partial G$, $G$ must contain the (unique) critical point $X^t$ of $F_t$, and $\partial G$ separates $X^t$ from $\partial J(\delta)$.
On the other hand, Lemma~\ref{LemmaG} together with Lemma~\ref{LemmaE} yield that
on one of the segments $I_i$, $i=1, 2$, the function $F_t$ does not change its sign.
The resulting contradiction proves the first part of the lemma.

To prove the second part, first, we notice that, since the sets $\{F_t>0\}\cap J(\delta)$ and $\{F_t<0\}\cap J(\delta)$
cannot be simultaneously connected, the set $\{F_t\ne 0\}\cap J(\delta)$ has at least three
connected components.
One of them, we call it $\Omega_0$, contains the critical point $X^t$ and therefore, by Lemma~\ref{LemmaG}, it contains one of the segments $I_i$. Since $F_t$ does not vanish on
$\partial^* J(\delta)$ (Lemma~\ref{LemmaE}), the boundary of $\Omega_0$
contains two opposite sides of $\partial^* J(\delta)$, the ones on which the end-points of the segment $I_i$ lie. For the same reason, there are two more connected components of
the set $\{F_t\ne 0\}\cap J(\delta)$, each of these two components
contains on its boundary one of two remaining opposite sides of the set
$\partial^* J(\delta)$. At last,  arguing as in the proof of the first part
(and using again Lemmas~\ref{LemmaF} and~\ref{LemmaG}), we see that the fourth
connected component of the set $\{F_t\ne 0\}\cap J(\delta)$ cannot exist.
\mbox{} \hfill $\Box$

\medskip\noindent\underline{Summary}:
Let $p$ be a saddle point of $f$. Suppose that conditions~\eqref{eq:parameters} and~\eqref{eq:delta} hold,
and let $J(p, \delta) = \Psi_p J(\delta)$ be the corresponding joint.
Suppose that $0$ is not a critical value of $f_t$.
\begin{itemize}
\item Then,
the set $Z(f_t)\cap J(p, \delta)$ consists of two connected components.
Each of them enters and exits the joint through its own terminals
$\Psi_p T_i$.
\item
We say that the joint $J(p, \delta)$ has {\em positive type} if
the set $J(p, \delta)\cap \{f_t>0\}$ is connected (and therefore,
the set $J(p, \delta)\cap \{f_t<0\}$ is disconnected and consists
of two connected components). Otherwise, we say that the joint $J(p, \delta)$
has {\em negative type}. Suppose that $|f_t(p)|\gtrsim A(\Delta\alpha)^2$.
Then, the type of the joint
$J(p, \delta)$ coincides with the sign of $f_t(p)$.
\end{itemize}

\subsection{Global matters: the gradient flow}
Fix the functions $f\in C^3(A, \Delta, \alpha, \beta)$ and $g\in C^3(A)$.
Let $f_t = f + t g$, $\widetilde{f}= f_{\alpha'}$,
and consider the gradient flow $z_t$, $0\le t \le \alpha'$, defined by the ODE
\begin{equation}\label{eq:ODE}
\frac{{\rm d}z_t}{{\rm d} t} = - \frac{(\partial_t f_t)(z_t)}{|\nabla f_t(z_t)|^2}\,
\nabla f_t(z_t)
\end{equation}
with the initial condition $ z_0 \in Z(f)$.

\begin{lemma}\label{LemmaK}
Suppose that $A\Delta^2\beta^2 \ll \alpha \ll (A\Delta)^{-2}\beta$
and $A\alpha' \ll \alpha$.
Let $\delta= c(A\Delta)^{-1}$ with sufficiently small constant $c>0$. Then, for any arc
$ I\subset Z(f)\setminus \bigcup_{p\in\Cr(\alpha)} D(p, 2\delta^2) $,
the flow $z_t$ provides a $C^1$-homotopy of $I$ onto an arc
$ \widetilde{I}\subset Z(\widetilde{f})\setminus \bigcup_{p\in\Cr(\alpha)} D(p, \delta^2) $.
Vice versa, for any arc
$ \widetilde{I}\subset Z(\widetilde{f})\setminus \bigcup_{p\in\Cr(\alpha)} D(p, 2\delta^2) $,
the inverse flow $z_{\alpha'-t}$ provides a $C^1$-homotopy of $\widetilde{I}$ onto
an arc $ I\subset Z(f)\setminus \bigcup_{p\in\Cr(\alpha)} D(p, \delta^2) $.
Moreover, these homotopies move the points by at most $O(A\alpha'/\beta)$.
\end{lemma}

\noindent{\em Proof of Lemma~\ref{LemmaK}}:
Let $X\stackrel{\rm def}=\bigl\{ |\nabla f| < \beta \bigr\}$ and
\[
\bar\Omega = \bar\Omega (\eps, \eta) \stackrel{\rm def}=
\bigl\{(t, x)\in [-\eps, \alpha'+\eps] \times
(\bS^2(L)\setminus X)\colon |f_t(x)|\le \eta \bigr\}
\]
(the choice of small positive parameters $\eps$ and $\eta$ has no importance), and let
$\Omega$ be the interior of $\bar\Omega$. Note that
\[
|\partial_t f_t | = |g| \le A \qquad {\rm everywhere},
\]
and
\[
|\nabla f_t | \ge |\nabla f| - A\alpha' \ge \tfrac12 \beta \qquad
{\rm on\ } \bS^2(L)\setminus X.
\]
Therefore, the flow moves the points with the speed
\[
\Bigr| \frac{{\rm d}z_t}{{\rm d}t} \Bigl| =
\Bigl| \frac{(\partial_t f_t)(z_t)}{|\nabla f_t(z_t)|}  \Bigr| \le \frac{2A}{\beta}.
\]

The RHS of the ODE~\eqref{eq:ODE} is a $C^1$-function on $\bar\Omega$.
Therefore, for any initial point $z_0\in Z(f)\setminus \bar X$, the ODE has a
unique $C^1$-solution. The solution  exists
until it reaches the boundary of $\Omega$. Note that along the
trajectory $z_t$ we have
\[
\frac{\rm d}{{\rm d}t}\, f_t(z_t) =
\bigl(\, \frac{\partial}{\partial t}  f_t\, \bigr)(z_t) +
\bigl\langle \nabla f_t(z_t), \frac{{\rm d}z_t}{{\rm d}t}\, \bigr\rangle = 0,
\]
whence $f_t(z_t)=0$ (recall that $z_0\in Z(f)$).
Hence, if the solution $z_t$ is not defined on $[0, \alpha']$, then there exists $\tau\le\alpha'$ such that $d_L(z_t, X)\to 0$ as $t\uparrow \tau$. This means that
the closure of the trajectory $z_t$, $0\le t < \tau$, contains a point $\bar z$ with
$|\nabla f(\bar z)|\le \beta$. Recalling that the point $z_t$ moves with the speed
$\le 2A/\beta$, we see that $ d_L(\bar z, z_0) \le 2A\alpha'/\beta$.
Furthermore, by the continuity of $f_t$, we have
$f_\tau(\bar z)=0$, whence  $|f(\bar z)|\le \alpha' |g(\bar z)| \le A\alpha' \ll \alpha$.
Combining this with the gradient estimate $|\nabla f(\bar z)| \le \beta$ and
applying Lemma~\ref{lemma:crit-almost_sing}, we conclude that there is a unique critical
point $p$ of $f$ with $d_L(p, \bar z) \le 2\beta\Delta$, i.e., with
\[
d_L(p, z_0) \le 2\beta\Delta + \frac{2A\alpha'}{\beta},
\]
which is much less than $\delta^2$.

It remains to check that this critical point $p$ belongs to $\Cr (\alpha)$, which is
straightforward:
\[
|f(p)| \lesssim |f(\bar z)| + d_L(\bar z, p)\beta + O(d_L(\bar z, p)^2 A)
\lesssim A\alpha' + \Delta\beta^2 + A \Delta^2\beta^2 \ll \alpha
\]
since $A\alpha'$ and $A\Delta^2\beta^2$ are both much less than $\alpha$.

Since the function $\widetilde{f} = f + \alpha' g$
belongs to the class $C^3(2A, 2\Delta, 2\alpha, 2\beta)$,
the same arguments can be also applied to the inverse flow $z_{\alpha'-t}$. \mbox{} \hfill $\Box$

\subsection{The upshot}
We start with functions
$f\in C^3(A, \Delta, \alpha, \beta)$ and $g\in C^3(A)$, and consider
the perturbation $\widetilde{f}=f+\alpha'g$.
We assume that the parameters
\[
\alpha'\ll \alpha\ll\beta\ll 1\ll A \ll \Delta
\]
satisfy the following relations:
\[
A\alpha' \ll \alpha, \quad A\Delta^2\beta^2 \ll \alpha \ll (A\Delta)^{-2}\beta,
\quad A^2 \Delta^3\alpha \ll 1
\]
(which, in particular, yield conditions~\eqref{eq:parameters}).
We also assume that
the perturbation $\widetilde{f}$ is not too small on $\Cr(\alpha)$:
\[
\min_{\Cr(\alpha)} |\widetilde{f}| \gtrsim A \Delta^2\alpha^2.
\]
We set
\begin{align*}
\operatorname{Cr}_{\tt S}(\alpha) &=
\{p\in\Cr(\alpha)\colon p\  {\rm is\ a\ saddle\ point\ of\ } f \}, \\
\operatorname{Cr}_{\tt E}(\alpha) &=
\{p\in\Cr(\alpha)\colon p\  {\rm is\ a\ local\ extremum\ of\ } f \}.
\end{align*}
We put $\delta= c(A\Delta)^{-1}$ with sufficiently small
positive constant $c$, and consider the disks $D(p, \delta)$,
$p\in\operatorname{Cr}_{\tt E}(\alpha)$,
and the joints $J(p, \delta)$, $p\in\operatorname{Cr}_{\tt S}(\alpha)$.
If the constant $c$ in the definition of $\delta$
was chosen sufficiently small, then all these disks and joints
are mutually disjoint (recall that by Lemma~\ref{LemmaB} the points from the
set $\Cr(\alpha)$ are $c_0 (A\Delta)^{-1}$-separated with a positive
constant $c_0$).

\subsubsection{Stable loops}

These are connected components of $Z(f)$ and $Z(\widetilde{f})$
that do not intersect the set
\[
U = U(\operatorname{Cr}(\alpha), \delta)
\stackrel{\rm def}=
\Bigl( \bigcup_{p\in \operatorname{Cr}_{\tt E}(\alpha)} D(p, \delta) \Bigr)\,
\bigcup\,
\Bigl( \bigcup_{p\in \operatorname{Cr}_{\tt S}(\alpha)} J(p, \delta) \Bigr).
\]
We denote by $N_{\rm I}(f)$ the number of stable loops in $Z(f)$ and by
$N_{\rm I}(\widetilde{f})$
the number of stable loops in $Z(\widetilde{f})$.

Observe that
$D\bigl( p, \sqrt{\tfrac{a}{b}}\, \delta \bigr) \subset J(p, \delta)$,
$p\in \operatorname{Cr}_{\tt S}(\alpha)$, and that
$\sqrt{\tfrac{a}{b}} \ge (A\Delta)^{-1/2}$, we see that, for each
$p\in \operatorname{Cr}(\alpha)$, we have $D(p, 2\delta^2)\subset U$.
Therefore, Lemma~\ref{LemmaK} applies to stable loops in $Z(f)$ as well as to stable loops
in $Z(\widetilde f)$ and yields a one-to-one correspondence between the set of
stable loops in $Z(f)$ and the set of stable loops in $Z(\widetilde{f})$.
That is, $N_{\rm I}(\widetilde{f})=N_{\rm I}(f)$.

\subsubsection{Blinking circles}

These are small connected components of $\widetilde{f}$
that surround the points $p\in \operatorname{Cr}_{\tt E}(\alpha)$ and lie in the interiors
of the corresponding disks $D(p, \delta)$. Recall that, by Lemma~\ref{LemmaC},
$Z(\widetilde f)$ cannot intersect the boundary circle
$\partial D(p, \delta)$ of such a disk.

By the summary in the end of the local extrema section~\ref{subsect:blinking-circles},
the number of such components is
\[
N_{\rm II}(\widetilde f) \stackrel{\rm def}=
\bigl| \bigl\{p\in\operatorname{Cr}_{\tt E}(\alpha)\colon
\widetilde{f}(p)\ {\rm and\ the\ eigenvalues\ of\ }
\cH_f(p)\ {\rm have\ opposite\ signs\,}  \bigr\} \bigr|.
\]

\subsubsection{The Bogomolny-Schmit loops}\label{subsubsect:BS-loops}

This is the most interesting part of $Z(\widetilde f)$.
Consider the graph $G=G(f)$ embedded in $\bS^2(L)$. The vertices of $G$ are
the joints $J(p, \delta)$, $p\in\operatorname{Cr}_{\tt S}(\alpha)$. The edges
are connected components of the set
\begin{equation}\label{eq:*}
Z(f) \setminus \bigcup_{p\in \operatorname{Cr}_{\tt S}(\alpha)} J(p, \delta)
\end{equation}
that touch the boundaries $\partial J(p, \delta)$ (these components are
homeomorphic to intervals, while the other connected components of the
set~\eqref{eq:*} are homeomorphic to circles). Each vertex of this graph has degree $4$.
The signs of $\widetilde{f}(p)$, $p\in \operatorname{Cr}_{\tt S}(\alpha)$,
determine the way the graph $G$ is turned into a collection of loops, see the summary in the
end of the saddle point section~\ref{subsect:saddle-points}.
\begin{figure}[h]
\includegraphics[width=0.5\textwidth]{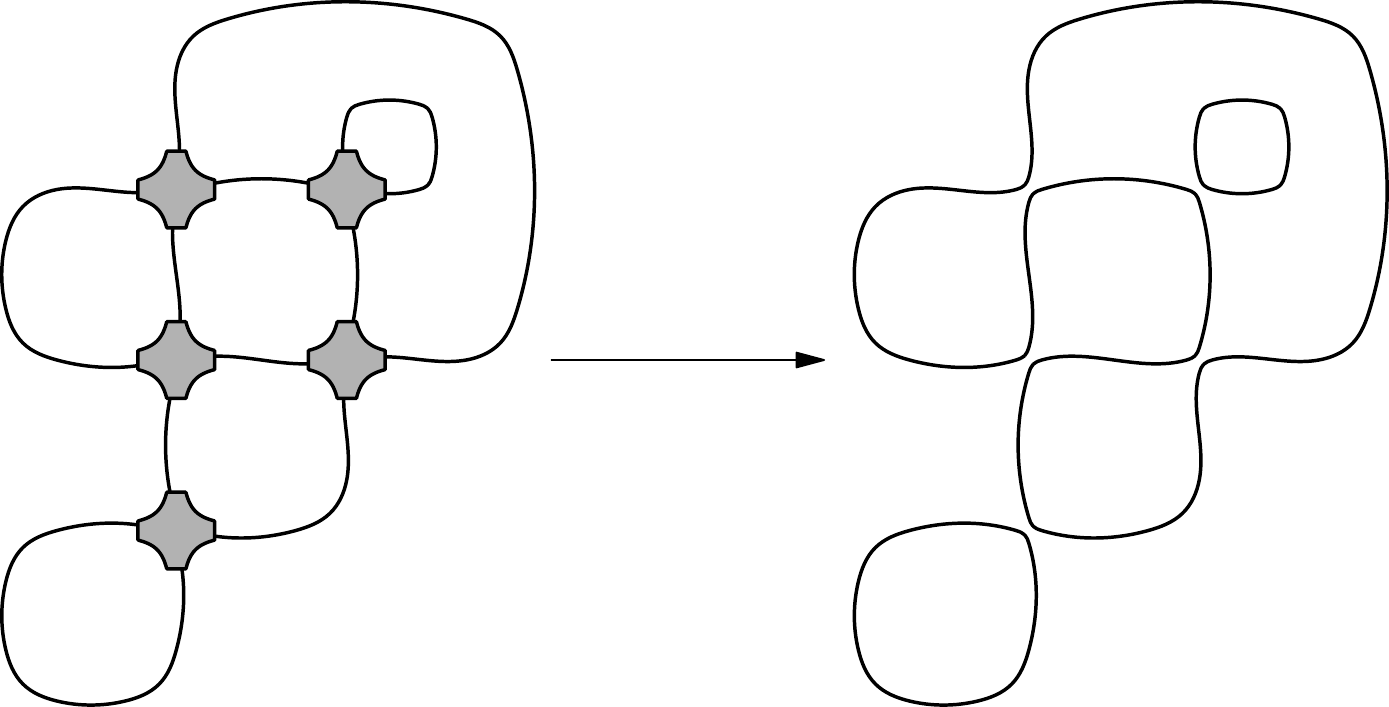}
\caption{Creation of the Bogomolny-Schmit loops}
\label{zfigure2}
\end{figure}
By $N_{\rm III}(\widetilde f)$ we denote the number of loops in this collection.

\subsubsection{}
At last, we are able to state the main result of this section:
\begin{lemma}\label{LemmaM}
Let $f\in C^3(A, \Delta, \alpha, \beta)$, $g\in C^3(A)$ and
$\widetilde{f}=f+\alpha'g$. Suppose that the parameters
\begin{equation}\label{eq:parameters0}
\alpha'\ll \alpha\ll\beta\ll 1\ll A \ll \Delta
\end{equation}
satisfy the following relations:
\begin{equation}\label{eq:parameters1}
A\alpha' \ll \alpha, \quad A\Delta^2\beta^2 \ll \alpha \ll (A\Delta)^{-2}\beta,
\quad A^2 \Delta^3\alpha \ll 1,
\end{equation}
and that
\begin{equation}\label{eq:parameters2}
\min_{\Cr(\alpha)} |\widetilde{f}| \gtrsim A \Delta^2\alpha^2.
\end{equation}
Then,
\[
N(\widetilde f\, ) = N_{\rm I}(\widetilde f\, ) + N_{\rm II}(\widetilde f\,) +
N_{\rm III}(\widetilde f\,).
\]
\end{lemma}

\section{Lower bounds for the Hessian of $f_L$ on the almost singular set}
Now, we return to regular Gaussian ensembles $(f_L)$.
%First, we show that, w.h.p., the random function $f_L$ belongs to the
%class $C^3(\log L, L^{3\e}, \alpha, \beta)$ with
%any $\alpha\le L^{-2+2\eps}$, $\beta \le L^{-1-\frac12 \eps}$.

\begin{lemma}\label{lemma-Hessian-lowerbound1}
Given a sufficiently small positive $\e$, let $\alpha\le L^{-2+2\e}$, $\beta^2 L^{3\e} \le \alpha$.
Then there exists $L_0=L_0(\e)$ such that, for each $L\ge L_0$, w.h.p.,
\[
\max_{\Cr(\alpha, \beta)} \bigl\| (\nabla^2 f_L)^{-1} \bigr\|_{\tt op} \le L^{3\e},
\]
where $\|\, . \,\|_{\tt op}$ denotes the operator norm.
\end{lemma}

\noindent{\em Proof}:
Fix a small $\eps>0$ and consider the set $\Cr(5\alpha, 4\beta)$.
Let $p$ be the probability that a given point $x\in\bS^2(L)$ belongs to the set
$\Cr(5\alpha, 4\beta)$. By the invariance of the ensemble
$(f_L)$, this probability does not depend on $x$.
The statistical independence of $f_L(x)$ and $\nabla f_L(x)$
(Lemma~\ref{lemma:independence}), non-degeneracy of their distributions
(Lemma~\ref{lemma:non-degeneracy}), and uniform boundedness of their
variances yield that
\[
p = \bP\bigl\{ |f_L(x)|\le 5\alpha \bigl\} \cdot \bP\bigl\{ |\nabla f_L(x)|\le 4\beta\bigr\}
\simeq \alpha \beta^2\,.
\]
Next, note that, by Fubini's theorem,
\[
\bE[\,\operatorname{area}( \Cr (5\alpha, 4\beta) )\,]
= p \cdot \operatorname{area}(\bS^2(L)) \simeq p\cdot L^2,
\]
whence, by Chebyshev's inequality,
\[
\bP\bigl\{\,\operatorname{area}( \Cr (5\alpha, 4\beta) ) > p\cdot L^2\cdot L^{\frac12 \e}\, \bigr\} \lesssim L^{-\frac12 \e}.
\]
Thus, w.h.p.,
\[
\operatorname{area}( \Cr (5\alpha, 4\beta) ) \le p\cdot L^2\cdot L^{\frac12 \e}
\simeq (\alpha L^2) \beta^2 L^{\frac12 \e} \stackrel{\alpha L^2 \le L^{2\e}}\le \beta^2 L^{2.5 \e}.
\]

Denote by $\mu = \mu (x)$ the eigenvalue of the Hessian matrix $\nabla^2 f_L$ with the
minimal absolute value, and by $w=w(x)$ the corresponding normalized eigenvector.
Assume that $\| f_L \|_{C^3}<\log L$ (recall that this holds w.o.p.)
and suppose that, for some $x\in\Cr (\alpha, \beta)$
and $L\ge L_0$, $ |\mu (x)|< \Delta^{-1} $, where $\Delta = L^{3\e}$.
We will show that then the set $\Cr (5\alpha, 4\beta)$ contains a subset
$\widetilde{G}=\widetilde{G}(x)$ with
$ \operatorname{area}(\widetilde{G}) \gtrsim \beta^2 L^{3\e} (\log L)^{-1}$. This
will immediately imply the lemma.

\medskip
Fix a point $x\in\Cr(\alpha, \beta)$ with $ |\mu (x)|< \Delta^{-1} $,
take the corresponding map $\Psi_x$, and let $F=f_L\circ \Psi_x$.
Put $\tau=\beta \Delta \ll 1$. For $Y=tw(x)$, $0\le t \le \tau$, we have
\begin{multline*}
|\nabla F(Y)| \ \le
|\nabla F(0)| + |\nabla^2 F(0) Y| + O(|Y|^2 \| F \|_{C^3(\{|X|\le 1\})}) \\
\le \beta + \Delta^{-1}\tau + \tau^2 O(\log L)
= 2\beta + \beta^2\Delta^2\, O(\log L) \le 3\beta.
\end{multline*}
Then, letting $I=[0, tw(x)]$, we get
\[
|F(Y)| \le
|F(0)| + \max_I |\nabla F| \cdot |Y|
\le \alpha + 3\beta \tau
= \alpha + 3 \beta^2\Delta \le  4\alpha\,,
\]
since $\beta^2 \Delta \le \alpha$.

Put $\rho=c \beta (\log L)^{-1}$ with a sufficiently small positive constant $c$,
and denote by $\Omega$ the $\rho$-neighbourhood of the segment $I$
on the plane $\Pi_x$. Then,
\[
\max_{\overline{\Omega}} |\nabla F| \le
\max_I | \nabla F | + O( \rho \| F \|_{C^3(\{|X|\le 1\})} )
\le 3\beta + \rho\, O(\log L) < 3.5\beta,
\]
and
\begin{multline*}
\max_{\overline\Omega} |F| \le
\max_I | F | + \rho \max_I | \nabla F | + O(\rho^2 \| F \|_{C^3(\{|X|\le 1\})} ) \\
\le 4\alpha + 3\beta \rho + \rho^2 O(\log L)
= 4\alpha +\beta^2\, o(1) < 5\alpha\,.
\end{multline*}
Let $\widetilde{\Omega}=\Psi_x(\Omega)$. We see that
$\widetilde{\Omega}\subset \Cr(5\alpha, 4\beta)$.

At the same time,
\[
\operatorname{area}(\widetilde{\Omega}) \gtrsim \tau \rho = c \beta^2 \Delta (\log L)^{-1}
= c \beta^2 L^{3\e} (\log L)^{-1},
\]
completing the proof.
\hfill $\Box$

\medskip
The next lemma gives us a lower bound for the probability that a given point
$x\in\bS^2(L)$ belongs to the set
\[
\Cr(\alpha, \beta, \Delta) = \{x\in\Cr(\alpha, \beta)\colon
\|(\nabla^2 f_L(x))^{-1}\|_{\tt op}\le\Delta\}.
\]
Let $p$ be the probability that a given point $x\in\bS^2(L)$ belongs to the set
$\Cr(\alpha, \beta)$. By the invariance of the ensemble $(f_L)$, this
probability does not depend on $x$.

\begin{lemma}\label{Hessian_conditional}
For any $\alpha, \beta \le 1$ and any $\Delta\ge 2$,
\[
\bP\bigl\{x\in\Cr(\alpha, \beta, \Delta)\bigr\}
\ge (1-C(\log \Delta)^{\frac14}\Delta^{-\frac12})p.
\]
\end{lemma}

\medskip\noindent{\em Proof}:
We need to show that, conditioned on $x\in\Cr(\alpha, \beta)$, the probability that
\[
\|(\nabla^2 f_L(x))^{-1}\|_{\tt op}\ge \Delta
\]
is $\lesssim (\log \Delta)^{\frac14}\Delta^{-\frac12}$.
Since $f_L(x)$ and its Hessian $\nabla^2 f_L(x)$ are independent of the
gradient $ \nabla f_L(x) $, it will suffice to show that
\[
\bP \bigl\{\|(\nabla^2 f_L(x))^{-1}\|_{\tt op} > \Delta\, \big|\, |f(x)|\le\alpha \bigr\}
\lesssim (\log\Delta)^{\frac14} \Delta^{-\frac12}.
\]

Fix $x\in\bS^2(L)$ and denote by
$\mu_1(x)$, $\mu_2(x)$ the eigenvalues of the Hessian. First, we show that,
conditioned on the event $\{ |f(x)|\le\alpha\}$,
with large probability, $|\mu_1 (x) \cdot \mu_2 (x)| = |\det \nabla^2 f_L(x)|$ cannot be
too small, and then that, with large probability,
$\max(|\mu_1(x)|, |\mu_2(x)|)= \| \nabla^2 f_L(x)\|_{\tt op} $ cannot be too big.
Together, these two estimates will do the job.

Fix coordinates $(X_1, X_2)$ in the plane $\Pi_x$.
Then
\[
\det \nabla^2 f_L(x) = \partial^2_{1, 1} f_L(x) \partial^2_{2, 2} f_L(x) - ( \partial^2_{1, 2} f_L(x) )^2.
\]
Recalling that, by Lemma~\ref{lemma:independence},
$ \partial^2_{1, 2} f_L(x) $ is independent of the vector
\[ \bigl( f_L(x), \partial^2_{1, 1} f_L(x), \partial^2_{2, 2} f_L(x) \bigr)^{\tt t}, \]
and that by Lemma~\ref{lemma:non-degeneracy}, the distribution of
$\partial^2_{1, 2} f_L(x) $ does not degenerate, we conclude that,
for any $\delta>0$,
\[
\bP\bigl\{ | \mu_1(x) \cdot \mu_2 (x)  | < \delta\, \big|\,
|f(x)|\le \alpha \bigr\}
\le \sup_{s\in\bR} \bP\bigl\{ | (\partial^2_{1, 2} f_L(x))^2-s | <  \delta \bigr\} \\
\lesssim \sqrt{\delta}.
\]

In the second step, taking into account that
$\|\nabla^2 f_L(x)\|_{\tt op}
\le 2\max_{1\le i, j \le 2} |\partial^2_{i, j} f_L(x) |$,
we need to estimate from above the absolute values of the second order
derivatives of $f_L$ at $x$ conditioned on $f_L(x)$.
The mixed derivative
$\partial^2_{1, 2} f_L(x)$ has a bounded variance and is independent of $f_L(x)$.
Furthermore, by the normal correlation theorem\footnote{
It says that if $(\theta, \xi)$ is a two-dimensional Gaussian vector,
then the expectation and variance of $\theta$ conditioned on $\xi$
equal
\[
\bE[\,\theta|\xi\,]= \frac{\bE[\,\theta \xi\,]}{\Var[\,\xi\,]}\, \xi
\]
and
\[
\Var[\,\theta|\xi\,] = \Var[\,\theta\,] - \frac{(\bE[\,\theta\xi\,])^2}{\Var[\,\xi\,]}
\le \Var[\,\theta\,].
\]
},
the distribution of
$\partial^2_{i, i} f_L(x)$, $i=1, 2$, conditioned on $f_L(x)$, is normal with
bounded conditional mean
\[
\bigl| \bE\bigl[ \partial^2_{i, i} f_L(x) \,|\, f_L(x) \bigr] \bigr| =
\bigl|\bE\bigl[f_L(x)\, \partial^2_{i, i} f_L(x)\bigr]\bigr|
\cdot \bigl |f_L(x) \bigr| \lesssim 1
\]
(recall that $f_L$ is normal and that we are interested only in the
values $|f_L(x)|\le \alpha \le 1$), and with the bounded conditional variance
\[
\Var\bigl[ \partial^2_{i, i} f_L(x) \,|\, f_L(x) \bigr]
\le \Var\bigl[ \partial^2_{i, i} f_L(x) \bigr] \lesssim 1.
\]
Thus, $ \bP\bigl\{ \max(|\mu_1(x)|, |\mu_2(x)|)>\lambda\, \big|\, |f(x)|\le \alpha \bigr\}
\lesssim e^{-c\lambda^2}$. Therefore, after conditioning on $\{ |f(x)|\le \alpha \}$,
with probability at least
$ 1 - C\bigl( \sqrt{\delta} + e^{-c\la^2} \bigr)$, we have
\begin{multline*}
\delta \le |\mu_1 (x)| \cdot |\mu_2(x)|
= \min(|\mu_1(x)|, |\mu_2(x)|) \cdot \max(|\mu_1(x)|, |\mu_2(x)|) \\
\le \la \min(|\mu_1(x)|, |\mu_2(x)|) = \la \| (\nabla^2 f_L(x))^{-1}\|^{-1}_{\tt op}.
\end{multline*}
Letting $\la = \delta \Delta $ and
$\delta = \tfrac1{\sqrt{2c}} \Delta^{-1}\sqrt{\log\Delta}$, and noting that
$\sqrt{\delta} + e^{-c\la^2} \lesssim \Delta^{-\frac12} (\log\Delta)^{\frac14}$,
we complete the proof. \hfill $\Box$.

\section{The two-point function of the set $\Cr(\alpha, \beta)$}

In this section, we will look
at the set $ \Cr (\alpha, \beta) $ of almost singular points of $f_L$:
\[
\Cr (\alpha, \beta) = \bigl\{x\in\bS^2(L)\colon |f_L(x)|\le \alpha, |\nabla f_L(x)|\le\beta \bigr\}
\]
with very small parameters $\alpha$ and $\beta$, and at its one- and two-point functions
$p$ and $p(x, y)$. As above, $p$ is the probability that a given point $x\in\bS^2(L)$ belongs to the set $\Cr(\alpha, \beta)$. Recall that, by the invariance of the ensemble
$(f_L)$, this probability does not depend on $x$, and that
the statistical independence of $f_L(x)$ and $\nabla f_L(x)$ yields that
\[
p = \bP\bigl\{ |f_L(x)|\le\alpha \bigl\}\, \cdot \, \bP\bigl\{ |\nabla f_L(x)|\le\beta\bigr\}
\simeq \alpha \beta^2\,.
\]
By $p(x, y)$ we denote the probability that two given points $x, y \in\bS^2(L)$ belong to the set $\Cr (\alpha, \beta)$. By the invariance of the distribution of $f_L$ with respect to isometries of the sphere, $p(x, y)$ depends only on the spherical distance between the points $x$ and $y$.

\subsection{Estimates of the two-point function}

\begin{lemma}\label{lemma:2point}
Let $(f_L)$ be a regular Gaussian ensemble. Then
the following estimates hold uniformly in $\alpha, \beta \le 1$
and in $L\ge L_0$:
\begin{equation}\label{eq:2point_short}
p(x, y) \lesssim \max\bigl\{d_L(x, y)^{-\Theta}, 1\bigr\} p^2
\end{equation}
with some positive constant $\Theta$, and for $d_L(x, y) \ge 1$, we have $p(x, y)=W p^2$, with
\begin{equation}\label{eq:2point_long}
|W-1| \lesssim d_L(x, y)^{-\gamma}.
\end{equation}

\end{lemma}

Note that the proof of the short-distance estimate~\eqref{eq:2point_short}
is quite lengthy, while the long-distance estimate~\eqref{eq:2point_long}
is a straightforward consequence of the power decay of correlations of $(f_L)$.

\subsection{Proof of the short-distance estimate~\eqref{eq:2point_short}}

\subsubsection{Beginning the proof}
Let $x, y\in\bS^2(L)$.
Fix the coordinate systems in the planes $\Pi_x$ and $\Pi_y$, and
let $\Gamma (x, y)$ be the covariance matrix of the Gaussian six-dimensional vector
\[
v(x, y) = \bigl(f_L(x), \nabla f_L (x), f_L(y), \nabla f_L(y) \bigr)^{\tt t}.
\]
Then
\[
p(x, y) =  \frac1{(2\pi)^3\,\sqrt{\det \Gamma (x, y)}}\,
\int_\Omega \exp\bigl[ -\tfrac12 \xi^{\tt t} \Gamma (x, y)^{-1} \xi\bigr]\, {\rm d}^6\xi,
\]
where
\[
\Omega = \bigl\{\xi\in\bR^6\colon
|\xi_1|, |\xi_4| \le\alpha, |\xi_2|^2+|\xi_3|^2, |\xi_5|^2+|\xi_6|^2 \le\beta^2
\bigr\}.
\]
Note that $\vol(\Omega)\lesssim \alpha^2 \beta^4 \simeq p^2 $. Hence,
to prove estimate~\eqref{eq:2point_short} we need to bound from below the minimal
eigenvalue $\la=\la (x, y)$ of the covariance matrix $\Gamma=\Gamma (x, y)$,
\[
\la = \min\bigl\{\xi^{\tt t} \Gamma (x, y) \xi\colon
\xi\in\bR^6, |\xi|=1 \bigr\}.
\]
Note that $\la$ does not depend on the choice of the coordinate systems
in the planes $\Pi_x$ and $\Pi_y$.

First, we show that there exists a sufficiently large constant $d_0$,
independent of $L$, so that $\la(x, y)$ is bounded from below by a positive
constant whenever $d_L(x, y)\ge d_0$. Hence, proving estimate~\eqref{eq:2point_short},
we assume that $d_L(x, y)\le d_0$
(while later, proving the long-distance estimate~\eqref{eq:2point_long},
we will assume that $d_L(x, y)\ge d_0$). The value of sufficiently large
constant $d_0$ is inessential for our purposes.

Denote by $\mathfrak C$ the covariance matrix of $\nabla f_L(x)$, and put
\[
\widetilde{\Gamma} =
\left(
\begin{matrix}
1 & 0 & 0 & 0 \\
0 & \mathfrak C & 0 & 0 \\
0 & 0 & 1 & 0 \\
0 & 0 & 0 & \mathfrak C
\end{matrix}
\right).
\]
Since the matrix $\mathfrak C$ is non-degenerate uniformly in $L\ge L_0$
(Lemma~\ref{lemma:non-degeneracy}),
the matrix $\widetilde\Gamma$ is also non-degenerate uniformly in $L\ge L_0$.

By our assumption on the power decay of correlations, we have
\[
\max_{1\le i,j \le 6}\ \bigl| \Gamma_{ij}(x, y) - \widetilde{\Gamma}_{ij} \bigr|
\lesssim (1 + d_L(x, y))^{-\gamma}.
\]
Then, provided that $d_L(x, y) \ge d_0$ with $d_0\gg 1$, we get
\[
\bigl\| \Gamma (x, y)^{-1} - \widetilde{\Gamma}^{-1}\bigr\|_{\tt op}
\lesssim d_L(x, y)^{-\gamma},
\]
and therefore,
\[
\bigl\| \Gamma (x, y)^{-1} \|_{\tt op} \ge \| \widetilde{\Gamma}^{-1}\bigr\|_{\tt op}
- O\bigl( d_L(x, y)^{-\gamma} \bigr) \ge \tfrac12\,
\| \widetilde{\Gamma}^{-1}\bigr\|_{\tt op}\,.
\]
Thus, till the end of the proof of the short-distance estimate~\eqref{eq:2point_short},
we assume that $d_L(x, y) \le d_0$ with some positive $d_0$ independent of $L$.

\medskip
Let $v(x)$ denote the three-dimensional  Gaussian vector
$ v(x) = \bigl( f_L(x), \nabla f_L(x) \bigr)^{\tt t} $,
and let $a=(\xi_1, \xi_2, \xi_3)^{\tt t}$, $b=(\xi_4, \xi_5, \xi_6)^{\tt t}$. Then
\[
\xi^{\tt t} \Gamma (x, y) \xi = \bE\bigl[ \langle v(x, y), \xi \rangle^2 \bigr]
= \bE\bigl[ ( \langle v(x), a \rangle  + \langle v(y), b \rangle)^2 \bigr],
\]
whence,
\[
\la = \min\bigl\{ \bE\bigl[ ( \langle v(x), a \rangle
+ \langle v(y), b \rangle)^2 \bigr]\colon
a, b\in\bR^3, |a|^2 +|b|^2=1 \bigr\}.
\]
By the compactness of the unit sphere in $\bR^3$,
there exist $a, b\in\bR^3$, $|a|^2+|b|^2=1$, such that
\begin{equation}\label{eq:lambda!}
\la = \bE\bigl[ ( \langle v(x), a \rangle  + \langle v(y), b \rangle)^2 \bigr]
= \| \langle v(x), a \rangle + \langle v(y), b \rangle \|^2.
\end{equation}
Here and till the end of the proof of Lemma~\ref{lemma:2point},
$\|\, .\,\|$ stands for the $L^2$-norm, i.e., $\|\eta\|^2=\bE[|\eta|^2]
=\Var[\eta]$ for the Gaussian random variable $\eta$.

\subsubsection{The normal and the tangential derivatives}
Let $\mathcal C$ be the big circle on $\bS^2(L)$ that
passes through the points $x$ and $y$, and let $\mathcal I\subset\mathcal C$ be
the shortest of the two arcs of $\mathcal C$ with the endpoints $x$ and $y$.
We orient $\mathcal C$ by moving from $x$ to $y$ along $\mathcal I$,
and choose the coordinate systems in $\Pi_x$ and $\Pi_y$ so that
one coordinate vector is parallel to the tangent to $\mathcal C$
(at $x$ and $y$ correspondingly), while the other one is orthogonal to $\mathcal C$.
We keep the same orientation for both coordinate systems.
We denote by $\partial_{\parallel}$ the derivative along $\mathcal C$ and
by $\partial_{\perp}$ the derivative in the normal direction to $\mathcal C$,
and decompose
\[
v(x) = v_{\parallel}(x) + v_{\perp}(x),
\]
where $v_{\parallel}(x) = \bigl( f_L(x), \partial_{\parallel} f_L(x), 0\bigr)^{\tt t}$
and $v_{\perp}(x) = \bigl(0, 0, \partial_{\perp} v(x) \bigr)^{\tt t}$.
Then, by~\eqref{eq:lambda!},
\begin{equation}\label{eq:lambda1}
\la = \| \langle v_{\parallel}(x), a' \rangle + \langle v_{\perp}(x), a'' \rangle +
\langle v_{\parallel}(y), b' \rangle + \langle v_{\perp}(y), b'' \rangle \|^2,
\end{equation}
where $a' = (a_1, a_2, 0)^{\tt t}$, $a'' = (0, 0, a_3)^{\tt t}$, similarly for $b'$ and $b''$,
and $|a'|^2 + |a''|^2 + |b'|^2 + |b''|^2 = 1$.

Now, consider another Gaussian six-dimensional
vector
\[
\widetilde{v}(x, y) = ( f_L(x),\, \partial_{\parallel} f_L(x), \, - \partial_\perp f_L(x),\,
f_L(y),\, \partial_{\parallel} f_L(y),\, - \partial_\perp f_L(y) )^{\tt t}.
\]
Since the distribution of $f_L$ is invariant with respect to orthogonal transformations,
the Gaussian vectors
$v(x, y)$ and $\widetilde{v}(x, y)$ have the same covariance matrix in the chosen
coordinate systems in $\Pi_x$ and $\Pi_y$. Therefore,
\begin{equation}\label{eq:lambda2}
\la = \| \langle v_{\parallel}(x), a' \rangle - \langle v_{\perp}(x), a'' \rangle +
\langle v_{\parallel}(y), b' \rangle  - \langle v_{\perp}(y), b'' \rangle \|^2.
\end{equation}
Juxtaposing \eqref{eq:lambda1} with \eqref{eq:lambda2}, we conclude that
\[
\la = \| \langle v_\parallel (x), a' \rangle + \langle v_\parallel (y), b' \rangle \|^2
+ \| \langle v_\perp (x), a'' \rangle + \langle v_\perp (y), b'' \rangle \|^2.
\]
We split the rest of the proof of estimate~\eqref{eq:2point_short} into two cases:
(i) $|a''|^2+|b''|^2\ge\tfrac12$ and (ii) $|a'|^2+|b'|^2\ge\tfrac12$.

\subsubsection{Case (i): $\large\bf|a''|^2+|b''|^2\ge\tfrac12$}

In this case, we use the estimate
$ \la \gtrsim \| \langle v_\perp (x), a'' \rangle + \langle v_\perp (y), b'' \rangle\|^2 $.
By the invariance of the distribution of $f_L$, the random pairs
$(\partial_\perp f_L(x), \partial_\perp f_L(y) )$ and
$(\partial_\perp f_L(y), \partial_\perp f_L(x) )$ have the same distribution.
Therefore,
\[
\| \langle v_\perp (x), a'' \rangle + \langle v_\perp (y), b'' \rangle \|
= \| \langle v_\perp (y), a'' \rangle + \langle v_\perp (x), b'' \rangle \|.
\]
Since $|a''|^2+|b''|^2\ge\tfrac12$, at least one of the following
holds:
\begin{itemize}
\item either $|a''+b''|\ge \tfrac12$, or $|a''-b''| \ge \tfrac12$.
\end{itemize}
We assume that $|a''+b''|\ge \tfrac12$ (the other case is similar and slightly
simpler), let $e=a''+b''$, $|e|\ge \tfrac12$, and notice that
\begin{multline*}
\| \langle v_\perp (x) + v_\perp (y), e \rangle \|
= \| \bigl( \langle v_\perp (x), a'' \rangle + \langle v_\perp (y), b'' \rangle \bigr)
+ \bigl( \langle v_\perp (x), b'' \rangle + \langle v_\perp (y), a'' \rangle \bigr) \| \\
\le
\| \langle v_\perp (x), a'' \rangle + \langle v_\perp (y), b'' \rangle \|
+ \| \langle v_\perp (x), b'' \rangle + \langle v_\perp (y), a'' \rangle \|
\lesssim \sqrt{\la}.
\end{multline*}
We take the point $z\in\mathcal C$, $z\ne x$, so that $d_L(y, z)=d_L(x, y)$. Then,
by the invariance of the distribution of $f_L$ with respect to the isometries of the sphere,
\[
\| \langle v_\perp (y) + v_\perp(z), e \rangle \|
= \| \langle v_\perp (x) + v_\perp(y), e \rangle \|,
\]
whence,
\[
\sqrt{\la} \gtrsim \| \langle v_\perp (x) - v_\perp(z),  e \rangle \|.
\]

Put $x_0=x$, $x_1=z$, then take the point $x_2\in \mathcal C$, $x_2\ne x_0$ so that
$d_L(x_2, x_1) = d_L(x_1, x_0)$, and continue this way till $d_L(x_0, x_N)\ge d_0$, where $d_0$
is the correlation length defined above. Then
\[
\| \langle v_\perp (x_0) - v_\perp(x_N), e \rangle\| \lesssim N \sqrt{\la}.
\]
On the other hand, since $d_0$ is the correlation length and $d_L(x_0, x_N)\ge d_0$,
we have
\[
\| \langle v_\perp (x_0) - v_\perp(x_N),  e \rangle \|^2
\gtrsim \| \langle v_\perp (x_0), e \rangle  \|^2 +
\| \langle v_\perp (x_N), e \rangle \|^2 \gtrsim 1
\]
(at the last step we use that the distribution of $v_\perp$ does not degenerate
uniformly in $L\ge L_0$ and that $|e|\ge \tfrac12$).
Therefore, $\la \gtrsim N^{-2}$.

Recalling that by the definition of $N$, we have $(N-1)d_L(x, y) < d_0$, we get
$\la \gtrsim d_L(x, y)^{2}$, which concludes our consideration of the first
case.

\subsubsection{Case (ii): $\large\bf |a'|^2+|b'|^2\ge\tfrac12$}
In this case, we restrict the function $f_L$ to the
big circle $\mathcal C$ and treat it as a periodic random Gaussian function
$F\colon \bR \to \bR$ with translation-invariant distribution. To simplify
the notation, we omit the index $L$. By $\rho$
we denote the spectral measure of $F$, that is,
\[
\bE [F(X)F(Y)] = \widehat{\rho}(X-Y)
= \int_{\bR} e^{2\pi {\rm i} \xi (X-Y)}\, {\rm d}\rho(\xi),
\]
where $X, Y\in\bR$ correspond to the points $x, y \in \mathcal C$. Then we have
\begin{align*}
\la &\gtrsim \bE\bigl[ \bigl( \langle v_{||}(X), a' \rangle
+ \langle v_{||}(Y), b' \rangle \bigr)^2 \bigr] \\
&= \int_{\bR} \bigl| (a_1 + a_2 \xi)e^{2\pi {\rm i}\xi X}
+ (b_1 + b_2 \xi)e^{2\pi {\rm i}\xi Y} \bigr|^2\, {\rm d}\rho(\xi) \\
&= \int_\bR \bigl|
(a_1 + a_2\xi)+(b_1 + b_2 \xi)e^{{\rm i}\delta \xi} \bigr|^2\, {\rm d}\rho(\xi)\,,
\end{align*}
where $\delta = \tfrac1{2\pi}\, |X - Y| = \tfrac1{2\pi}\, d_L(x, y)$. Furthermore,
\[
\rho (\bR) =\bE [F(0)^2] =1,
\]
\[
m \stackrel{\rm def}= \int_{\bR} \xi^2\, {\rm d}\rho(\xi) = E[F'(0)^2],
\quad m \simeq 1
\]
by the uniform non-degeneracy of $\nabla f_L$, and
\[
|\widehat{\rho}(s)| = |\bE[F(0)F(s)]| \lesssim |s|^{-\gamma}, \quad |s|\le \tfrac12 L,
\]
by the power decay of correlations of $f_L$. Note that since the function $F$ is
$2\pi L$-periodic, the Fourier transform $\widehat{\rho}$ of its spectral measure
is also $2\pi L$-periodic.

Recalling that
$\tfrac12 \le |a'|^2 + |b'|^2 \le 1$, we notice that
$|a_1|+|a_2|+|b_1|+|b_2| \ge \tfrac12 $ as well.
These remarks reduce the lower bound for $\la$ we are after
to a question in harmonic analysis.

Given $\delta>0$, consider the exponential sum\footnote{
Recall that the exponential sum is the expression
\[
S(\xi) = \sum_{j=1}^n q_j(\xi) e^{{\rm i}\la_j\xi},
\]
where $\la_j$ are real numbers and $q_j$ are polynomials in $\xi$ with complex coefficients.
As usual, $\deg S \stackrel{\rm def}= \sum_{j=1}^n (\deg q_j +1)$.
}
of degree $4$
\begin{equation}\label{eq:exp_sum}
p_\delta (\xi) = (a_1+a_2\xi) + (b_1+b_2\xi)e^{{\rm i}\delta\xi},
\qquad a_1, a_2, b_1, b_2\in\bC,
\end{equation}
and denote $\|p_\delta\|_W = |a_1|+|a_2|+|b_1|+|b_2|$. Then, the following lemma does the job.

\begin{lemma}\label{lemma:exp_poly}
Let $p_\delta$ be the exponential sum~\eqref{eq:exp_sum} of degree $4$.
Let $\rho$ be a probability measure on $\bR$ with the $2\pi L$-periodic Fourier transform $\widehat{\rho}$. Assume that
\[
m=\int_\bR \xi^2\, {\rm d}\rho (\xi) \simeq 1,
\]
and
\[
| \widehat{\rho}(s) | \le C|s|^{-\gamma}\,, \qquad |s|\le \tfrac12 L.
\]
Then, given $\delta_0>0$, there exist $c=c(C, \gamma, m, \delta_0)>0$ and
$L_0=L_0(C, \gamma, m, \delta_0)$ such that, for every $0<\delta\le \delta_0$ and every
$L\ge L_0$,
\[
\int_\bR
|p_\delta|^2\, {\rm d}\rho \ge c\delta^6 \| p_\delta \|_W^2.
\]
\end{lemma}

\subsection{Proof of Lemma~\ref{lemma:exp_poly}}

\subsubsection{Beginning the proof of Lemma~\ref{lemma:exp_poly}}
Let $A=2\sqrt m$. Then, by Chebyshev's inequality,
\[
\rho(\bR\setminus [-A, A]) \le A^{-2} m = \tfrac14\,.
\]
Take $B=\delta_0A$ and let $0<\delta \le \delta_0$. Then
$ [-A, A]\subset [-B\delta^{-1}, B\delta^{-1}]$.
Let $\kappa>0$ be a sufficiently small parameter, which we will choose later,
and consider the set
\[
\Xi = \bigl\{\xi\in [-B\delta^{-1}, B\delta^{-1}]\colon
|p_\delta(\xi)|<\|p_\delta\|_W\, \kappa^3 \delta^3 \bigr\}.
\]
We claim that
\begin{itemize}
\item[(a)] The set $\Xi$ is a union of at most $B+5$ intervals
$ I_j$, $1\le j \le B+5$.
\item[(b)] The length of each interval $I_j$ is
$\le C(B, \delta_0)\,\kappa$.
\end{itemize}
Having these claims, we will show that
$\rho\bigl( [-A, A]\setminus\Xi \bigr) \ge\tfrac14$, whence,
\[
\int_{[-A, A]\setminus\Xi } |p_\delta|^2\, {\rm d}\rho \ge \frac14 (\kappa\delta)^6
\| p_\delta \|^2_W.
\]
This will complete the proof of Lemma~\ref{lemma:exp_poly}.

\subsubsection{Proof of Claim (a)}
To show (a), we consider the exponential sum $P=|p_\delta|^2$
of degree $9$. By the classical Langer lemma
(see, for instance,~\cite[Lemma~1.3]{Nazarov}), the number of zeroes of any exponential sum
of degree $N$ on any interval $J\subset\bR$ cannot exceed
\[ (N-1) + \frac{\Delta}{2\pi}\, |J| \]
where $\Delta$ is the maximal distance between
the exponents in the exponential sum. Hence, the number of solutions
to the equation $P(\xi)=t$ on the interval $[-B\delta^{-1}, B\delta^{-1}]$ does not exceed
\[ 8+\frac{4\delta}{2\pi}\cdot 2B\delta^{-1}<8+2B. \]
Hence, the set $\Xi$ consists of at most $\tfrac12((8+2B)+2)= 5+B$ intervals, proving (a).

\subsubsection{Proof of Claim (b)}
To show (b), we apply Tur\'an's lemma~\cite[Theorem~1.5]{Nazarov},
which states that for any exponential sum $S$ of degree $N$ and any pair of closed
intervals $I\subset J$,
\[
\max_J |S| \le \left( \frac{C|J|}{|I|} \right)^{N-1}\, \max_I |S|\,,
\]
where $C$ is a numerical constant. Applying this lemma to the
exponential sum $p_\delta$
of degree $4$ and to each of
the intervals $I_j\subset [-B\delta^{-1}, B\delta^{-1}]$, we get
\begin{align*}
\max_{[-B\delta^{-1}, B\delta^{-1}]}\ |p_\delta| &\le
\left( \frac{C\cdot 2B\delta^{-1}}{|I_j|} \right)^3\, \max_{I_j} |p_\delta| \\
&\le \left( \frac{C\cdot 2B\delta^{-1}}{|I_j|} \right)^3\,
\|p_\delta\|_W\, \kappa^3\delta^3 \qquad\qquad  ({\rm since\ } I_j\subset\Xi) \\
&= \kappa^3 \left( \frac{C\cdot 2B}{|I_j|} \right)^3\, \|p_\delta\|_W.
\end{align*}
On the other hand,
\begin{align*}
\max_{[-B\delta^{-1}, B\delta^{-1}]}\ |p_\delta| &=
\max_{\xi\in [-B, B]}\, \bigl| (a_1+\delta^{-1}a_2\xi)
+ (b_1+\delta^{-1}b_2\xi) e^{{\rm i}\xi} \bigr| \qquad ({\rm by\ scaling}) \\
&\ge c(B) \bigl( |a_1|+\delta^{-1}|a_2|+ |b_1|+\delta^{-1}|b_2| \bigr)
\qquad ({\rm by\ compactness}) \\
&\stackrel{\delta\le \delta_0}\ge c(B) (1+\delta_0)^{-1} \| p_\delta \|_W.
\end{align*}
Thus,
$ |I_j| \le C(B, \delta_0) \kappa$,
proving (b).

\subsubsection{Completing the proof of Lemma~\ref{lemma:exp_poly}}
Recall that
$ \rho([-A, A]) \ge \tfrac34$, and that given $\kappa>0$, we defined
the set
\[
\Xi = \bigl\{\xi\in [-B\delta^{-1}, B\delta^{-1}]\colon
|p_\delta(\xi)|<\|p_\delta\|_W\, \kappa^3 \delta^3 \bigr\}
\]
satisfying (a) and (b).
Then we have the following alternative:
\begin{itemize}
\item
either $\rho([-A, A]\setminus \Xi)\ge \frac14$, or
$\rho([-A, A]\cap \Xi) \ge \frac 12$.
\end{itemize}
In the first case,
\[
\int_\bR |p_\delta|^2\, {\rm d}\rho
\ge \int_{[-A, A]\setminus\Xi} |p_\delta|^2\, {\rm d} \rho
\ge \tfrac14 \| p_\delta \|_W^2\, \kappa^6\delta^6.
\]
and we are done (modulo the choice of the parameter $\kappa$, which will be made
later). {\em It remains to show that if $\kappa$ is sufficiently small,
then the second case cannot occur}.

\medskip Suppose that $\rho([-A, A]\cap \Xi) \ge \frac 12$.
Then, $\rho(\Xi)\ge \frac12$. By claim (a), $\Xi$ is a union of at most $B+5$ intervals
$I_j$, hence, for at least one of them, $\rho (I_j) \ge c(B)>0$.
We call this interval $I$ and denote by $\nu$ the restriction of the measure $\rho$
on $I$. We choose a large parameter $S$ so that $1 \ll S \ll \kappa^{-1}$,
and estimate the integral
\[
J = \int_{-S}^S \Bigl( 1-\frac{|s|}S \Bigr) \bigl| \widehat{\nu}(s) \bigr|^2\,
{\rm d}s
\]
from below and from above, obtaining the estimates which will contradict each other.

Denote by $\xi_I$ the center of the interval $I$. Then,
for $|s|\le S$ and $\xi\in I$, we have
\begin{align*}
\bigl| e^{2\pi {\rm i}\xi s} - e^{2\pi {\rm i} \xi_I s }\bigr| &\le
2\pi |s| \cdot |\xi-\xi_I| \\
&\le 2\pi S \cdot \tfrac12 |I| \\
&\le 2\pi S \cdot \kappa \cdot \tfrac12 C(B, \delta_0)
\qquad ({\rm since,\ by\ claim\ (b)}, |I|\le C(B, \delta_0)\kappa ) \\
&< \tfrac12
\end{align*}
provided that $S\cdot\kappa$ is sufficiently small. Therefore, for $|s|\le S$, we have
\[
\bigl| \widehat{\nu} (s) \bigr|^2 = \left| \int_I e^{{\rm i}s\xi}\,
{\rm d}\rho(\xi) \right|^2
\ge \tfrac14 |\rho(I)|^2 \ge c_1(B)\,,
\]
whence,
\[
J \ge c_1(B) \int_{-S}^S \Bigl( 1 - \frac{|s|}{S} \Bigr) {\rm d}s  \ge c_2(B) S\,.
\]
On the other hand, using the identity
\[
\int_{\mathbb R} \phi (s) |\widehat{\nu}(s)|^2\, {\rm d}s
= \iint_{\bR\times\bR} \widehat{\phi}(\la-\eta)\, {\rm d}\nu(\la){\rm d}\nu(\eta),
\]
with $\phi (s) = (1-|s|/S)_+$ and noting that the Fourier transform of this function is
non-negative, we get
\[
\int_{-S}^S \Bigl( 1-\frac{|s|}S \Bigr) \bigl| \widehat{\nu}(s) \bigr|^2\, {\rm d}s
\le \int_{-S}^S \Bigl( 1-\frac{|s|}S \Bigr) \bigl| \widehat{\rho}(s) \bigr|^2\, {\rm d}s.
\]
Then, recalling that the function $\widehat{\rho}$ is $2\pi L$-periodic and that
$| \widehat{\rho}(s) | \le \min( 1, C|s|^{-\gamma})$ for $|s|\le L/2$,
and assuming without loss of generality that $\gamma<\frac12$, we get
\[
J \lesssim S\cdot (\min(S, L))^{-2\gamma}.
\]
Choosing $\kappa$ sufficiently small and $S$ sufficiently large, we arrive at a
contradiction, which completes the proof of Lemma~\ref{lemma:exp_poly}, and therefore, of
estimate~\eqref{eq:2point_short} in Lemma~\ref{lemma:2point}. \hfill $\Box$

\subsection{Proof of the long-distance estimate~\eqref{eq:2point_long}}

As above, we denote by $\mathfrak C$ the covariance matrix of $\nabla f_L(x)$, and put
\[
\widetilde{\Gamma} =
\left(
\begin{matrix}
1 & 0 & 0 & 0 \\
0 & \mathfrak C & 0 & 0 \\
0 & 0 & 1 & 0 \\
0 & 0 & 0 & \mathfrak C
\end{matrix}
\right).
\]
Since the matrix $\mathfrak C$ is non-degenerate uniformly in $L\ge L_0$,
the matrix $\widetilde\Gamma$ is also non-degenerate uniformly in $L\ge L_0$.

We assume that $d_L(x, y)\ge d_0$, where $d_0$ is sufficiently
large (and independent of $L$). Then we have
\[
\max_{1\le i,j \le 6}\ \bigl| \Gamma_{ij}(x, y) - \widetilde{\Gamma}_{ij} \bigr|
\lesssim d_L(x, y)^{-\gamma}.
\]
Therefore,
\[
\det \Gamma (x, y) =  \det\widetilde\Gamma  + O\bigl( d_L(x, y)^{-\gamma} \bigr)
\]
and
\[
\bigl\| \Gamma (x, y)^{-1} - \widetilde{\Gamma}^{-1}\bigr\|_{\tt op}
\lesssim d_L(x, y)^{-\gamma}.
\]

Recall that
\[
p(x, y) = \frac1{(2\pi)^{3} \sqrt{\det \Gamma (x, y)}}\,
\int_\Omega \exp\bigl[ -\tfrac12 \xi^{\tt t} \Gamma (x, y)^{-1} \xi\bigr]\, {\rm d}^6\xi,
\]
where
\[
\Omega = \bigl\{\xi\in\bR^6\colon
|\xi_1|, |\xi_4| \le\alpha, |\xi_2|^2+|\xi_3|^2, |\xi_5|^2+|\xi_6|^2 \le\beta^2
\bigr\},
\]
and that
\[
\frac1{(2\pi)^{3} \sqrt{\det\widetilde\Gamma}}\,
\int_\Omega
\exp\bigl[-\tfrac12 \xi^{\tt t} \widetilde{\Gamma}^{-1} \xi \bigr]\,
{\rm d}^6\xi = p^2.
\]
Hence,
\begin{align*}
p(x, y) &=
(2\pi)^{-3} (\det\widetilde\Gamma + O(d_L(x, y)^{-\gamma})^{-\frac12}\,
\int_\Omega
\exp\bigl[-\tfrac12 \xi^{\tt t} \widetilde{\Gamma}^{-1} \xi
+ O(d_L(x, y)^{-\gamma}) \bigr]\, {\rm d}^6\xi \\
&= (2\pi)^{-3}(1+O(d_L(x, y)^{-\gamma}) \, \frac1{\sqrt{\det\widetilde{\Gamma}}}\,
\int_\Omega
\exp\bigl[-\tfrac12 \xi^{\tt t} \widetilde{\Gamma}^{-1} \xi \bigr]\,
{\rm d}^6\xi \\
&= (1+O(d_L(x, y)^{-\gamma})\, p^2\,,
\end{align*}
completing the proof of estimate~\eqref{eq:2point_long} in Lemma~\ref{lemma:2point}.
\hfill $\Box$

\section{Structure of the set $\Cr(\alpha)$ with $L^{-2+\e}\le\alpha\le L^{-2+2\e}$}

Now, we are  ready to prove our main lemma:
\begin{lemma}\label{lemma_Cr(alpha)}
There exist positive $\eps_0$ and $c$, and positive $C$ such that, given $0<\eps\le \eps_0$ and $L\ge L_0(\eps)$, for
$L^{-2+\eps}\le\alpha\le L^{-2+2\eps}$, w.h.p.,
the set $\Cr(\alpha)$ is $L^{1-C\eps}$-separated, and $|\Cr(\alpha)|\ge L^{c\eps}$.
\end{lemma}

\subsection{W.h.p., the set $\Cr(\alpha)$ is $L^{1-C\e}$-separated}
In this part, we assume that $\alpha\le L^{-2+2\e}$, choose $\beta$ and $\rho$ so that
\[
\beta^2 L^{3\e} \le \alpha, \quad \rho=\beta(\log L)^{-1},
\]
and fix a maximal $\rho$-separated set $\mathcal X(\rho)$ on $\bS^2(L)$. Then
$|\mathcal X(\rho)|\simeq (L/\rho)^2$.

\medskip
First, we note that, w.h.p., the points
of the set $\Cr (\alpha)$ are $L^{-4\e}$-separated. This is a straightforward consequence
of part (B) in Lemma~\ref{lemma:crit-almost_sing} combined with a priori w.o.p.-bound
$\| f_L \|_{C^3} < \log L$ and with the w.h.p.-estimate
$ \max_{\Cr (\alpha)} \| (\nabla^2 f_L)^{-1} \|_{\tt op} \le L^{3\e} $
provided by Lemma~\ref{lemma-Hessian-lowerbound1}.
Hence, we need to estimate the probability of the event
\[
\mathcal E = \bigl\{ \exists z_1, z_2\in \Cr(\alpha)\colon L^{-4\e}
\le d_L(z_1, z_2) \le L^{1-C\e} \bigr\}
\]
with an appropriately chosen constant $C$.

Suppose that the event $\mathcal E$ occurs.
Denote by $x_1, x_2$ the closest to $z_1$, $z_2$  points in $\mathcal X(\rho)$.
Then,
\begin{equation}\label{eq:star}
\frac12 L^{-4\e} \le d_L(x_1, x_2) \le 2 L^{1-C\e}\,,
\end{equation}
and
\begin{gather*}
|f_L(x_i)| \le \alpha + O(\rho^2 \| f_L \|_{C^2}) \stackrel{\rm w.o.p.}<
\alpha + O(\beta^2 (\log L)^{-1}) < 2\alpha\,, \\
|\nabla f_L(x_i)| \le O(\rho \| f_L \|_{C^2}) \stackrel{\rm w.o.p.}<
\beta (\log L)^{-1} \cdot \log L = \beta\,,
\end{gather*}
i.e., $x_1, x_2 \in\Cr(2\alpha, \beta)$.
We claim that
\begin{itemize}
\item {\em the mean number of pairs of points
$x_1, x_2\in \Cr(2\alpha, \beta)\bigcap \mathcal X(\rho)$ satisfying~\eqref{eq:star} is
bounded from above by} $L^{-\e}$.
\end{itemize}
By Chebyshev's inequality, this yields that the
probability that there exists at least one such pair is also bounded from above by
$L^{-\e}$, which proves the $L^{1-C\e}$-separation.

\medskip
The mean we need to estimate equals
\[
\sum_{\substack{ x_1, x_2 \in\Cr(2\alpha, \beta)\cap\mathcal X(\rho) \\
\eqref{eq:star} {\rm\ occurs} }}\ p(x_1, x_2)\,,
\]
where $p(x_1, x_2)=\bP\{ x_1, x_2\in\Cr(2\alpha, \beta) \}$ is the
two-point function estimated in Lemma~\ref{lemma:2point}.
By Lemma~\ref{lemma:2point}, $p(x_1, x_2) \lesssim L^{4\e\Theta} p^2$, so
the whole sum is
\begin{align*}
&\lesssim \sum_{x_1\in\mathcal X(\rho)} L^{4\eps\Theta}p^2 \cdot
\bigl| \{x_2\in\mathcal X(r)\colon d_L(x_1, x_2)\le 2L^{1-C\epsilon}\}\bigr|
\\
&\lesssim L^{4\eps\Theta}p^2 \cdot (L^{1-C\eps}/\rho)^2 \cdot (L/\rho)^2 \\
%&\lesssim L^{-(C-4\Theta)\epsilon} p^2 \cdot (L/\rho)^4 \\
&\lesssim L^{-(C-4\Theta)\eps} \alpha^2 \beta^4 \cdot (L^4/\beta^4) \log^4 L
\qquad \qquad \qquad \qquad \qquad (p\simeq \alpha\beta^2, \rho = \beta/\log L) \\
&\lesssim L^{-(C-4\Theta)\eps} \cdot L^{4\eps}\log^4 L
\qquad \qquad \qquad \qquad \qquad \qquad \qquad
(\alpha \le L^{-2+2\eps}) \\
&< L^{-\eps},
\end{align*}
provided that the constant $C$ is sufficiently large.
This proves that the set $\Cr (\alpha)$ is $L^{1-C\e}$-separated.
\hfill $\Box$

\subsection{W.h.p., $|\Cr(\alpha)|\ge L^{c\e}$}

In this part we assume that $\alpha\ge L^{-2+\e}$ and introduce the parameters
$\beta$, $\Delta$ and $r$ satisfying
\[
\Delta=L^{\frac14 \e}, \quad \beta^2 L^{3\e} = \tfrac13 \alpha,
\quad r=\beta\Delta.
\]
We fix a maximal $r$-separated set $\mathcal X(r)$ on $\bS^2(L)$. Then
the disks $D(x, r)$, $x\in \mathcal X(r)$, cover $\bS^2(L)$ with a bounded multiplicity of
covering, and $|\mathcal X(r)|\simeq (L/r)^2$. We set
\[
Y = \mathcal X(r) \cap \Cr(\tfrac13 \alpha, \beta, \tfrac12\Delta).
\]
W.o.p., for $L\ge L_0$, we have $\| f_L \|_{C^3} \le \log L$. Then,
by Lemma~\ref{lemma:crit-almost_sing}, each disk $D(y, r)$, $y\in Y$, contains a unique critical point $z\in\Cr(\alpha)$, and
\begin{multline*}
\|(\nabla^2 f_L(z))^{-1}\|_{\tt op}\le \|(\nabla^2 f_L(y))^{-1}\|_{\tt op}
\bigl( 1 + O(r\, \|(\nabla^2 f_L(y))^{-1}\|_{\tt op} \cdot \| f_L \|_{C^3}) \\
\le \tfrac12\, \Delta + O(r\, \Delta^2 \log L) < \Delta,
\end{multline*}
provided that $L\ge L_0$.
This yields two useful observations which hold w.o.p.:
\begin{itemize}
\item[(i)] $|\Cr(\alpha)|\gtrsim |Y|$;
\item[(ii)] if $y_1, y_2\in Y$, then either
$d_L(y_1, y_2)\le 2r$ and the number of such pairs $(y_1, y_2)$
is $\lesssim|\mathcal X(r)|$, or $d(y_1, y_2)\ge L^{-\e}$.
Indeed, if the points $y_1, y_2\in Y$ generate the same critical point
$z$, then $d_L(y_1, y_2)\le 2r$. If they generate different critical points $z$,
we note that, by part~(B) of Lemma~\ref{lemma:crit-almost_sing},
the set of critical points $z$ of $f_L$ with $\|(\nabla^2 f_L(z))^{-1}\|_{\tt op}\le\Delta$ is
$c\Delta^{-1} (\log L)^{-1}$-separated, thus, in this case $d(y_1, y_2)\ge L^{-\e}$.
\end{itemize}
%In the $2$nd observation we used that, w.o.p., the set of critical points $z$
%of $f_L$ with $\|(\nabla^2 f_L(z))^{-1}\|_{\tt op}\le\Delta$ is
%$c\Delta^{-1} (\log L)^{-1}$-separated.

In what follows, we will show that,
for sufficiently large $L$,
\[
\bE[\, |Y| \,]\gtrsim L^{\frac12 \e}
\]
and that
\[
\Var[\, |Y|\, ] \lesssim L^{-c\e} (\bE [|Y|])^2.
\]
These two estimates combined with  the first observation readily yield what we need.

\subsubsection{Estimating $\bE[\, |Y| \,]$}
This estimate is straightforward:
\begin{align}
\nonumber
\bE[\, |Y| \,] &= \bP\bigl\{\, x\in\Cr(\tfrac13\alpha, \beta, \tfrac12\Delta) \,\bigr\}
\cdot |\mathcal X(r)| \\
\nonumber
&\gtrsim \bP\bigl\{\, x\in\Cr(\tfrac13\alpha, \beta) \,\bigr\}
\cdot \Bigl( \frac{L}r \Bigr)^2
\qquad \qquad ({\rm by\ Lemma~\ref{Hessian_conditional}})\\
\label{eq:E}
&\gtrsim \alpha\beta^2 \cdot \frac{L^2}{\beta^2\Delta^2} \\
\label{eq:E1}
&\ge L^{\frac12 \e} \qquad\qquad\qquad\qquad\qquad\qquad\quad
(\alpha\ge L^{-2+\e}, \Delta=L^{\frac14 \e}).
\end{align}

\subsubsection{Estimating $\Var[\, |Y| \,]$}

In this section, $p$ and $p(x, y)$ will denote the one- and two-point functions of the set
$\Cr\bigl( \tfrac13\alpha, \beta \Bigr)$. Given $x, y\in\bS^2(L)$, we put
\[
p_\Delta = \bP\bigl\{ y\in\Cr\bigl(\tfrac13\alpha, \beta, \tfrac12\Delta\bigr) \bigr\},
\quad p_\Delta (x, y) =
\bP\bigr\{x, y\in\Cr\bigl(\tfrac13\alpha, \beta, \tfrac12\Delta\bigr) \bigr\}.
\]
Then $ \bE[|Y|]=p_\Delta |\mathcal X(r)| $ and
\[
\Var[\, |Y|\,] = \sum_{y\in\mathcal X(r)} (p_\Delta - p_\Delta^2)
+ \sum_{\substack{x, y\in \mathcal X(r)\\ x\ne y }} (p_\Delta (x, y) - p_\Delta^2).
\]
The first sum on the RHS is bounded by
\[
p_\Delta\, |\mathcal X(r)| = \bE[\,|Y|\,] \stackrel{\eqref{eq:E1}}\lesssim
L^{-\frac12 \e}\bigl( \bE[|Y|]\bigr)^2.
\]
So we need to estimate the double sum only.

In the double sum we consider separately the terms with $d_L(x,y)\le 2r$, the terms with
$2r < d_L(x, y) < L^{-\e}$, the terms with $L^{-\e}\le d_L(x,y)\le L^{\e}$, and the terms
with $d_L(x, y) \ge L^\e$.

\medskip
\paragraph{\underline{The terms with $d_L(x,y)\le 2r$}.}
Taking into account that the number of such
pairs is $\lesssim |\mathcal X(r)|$, we bound this sum by
$ \lesssim p_\Delta |\mathcal X(r)| = \bE[|Y|] \stackrel{\eqref{eq:E1}}\lesssim
L^{-\frac12 \e}(\bE[|Y|])^2 $.

\medskip
\paragraph{\underline{The terms with $2r< d_L(x,y) <L^{-\e}$}.}
By the second observation, we conclude that w.o.p. this case
cannot occur, that is, the probability that there exists a pair of almost-singular points
$x, y\in\Cr(\tfrac13\alpha, \beta, \tfrac12 \Delta)$
with  $2r< d_L(x,y) <L^{-\e}$ is $O(L^{-C})$ with any positive $C$.
Thus,
in this range, $p_\Delta (x, y) = O(L^{-C})$, while the total number of pairs
$x, y$ is bounded by $|\mathcal X(r)|^2 \simeq (L/r)^4 \ll (L/\beta)^4 \ll L^{10}$,
provided that $\e$ in the definition of the parameters $\beta$ and $\alpha$ is sufficiently
small. That is, the sum is negligibly small.

\medskip
\paragraph{\underline{The terms with $L^{-\e} \le d_L(x, y)\le L^{\e}$}.}
In this case, we estimate each summand in the double sum by $p(x, y)$, which, by the short-distance estimate in Lemma~\ref{lemma:2point}, is
$\lesssim \max(d_L(x, y)^{-\Theta}, 1) p^2 \lesssim L^{\Theta\e} p_\Delta^2$.
Then, the whole double sum is
\[
\ll L^{\Theta\e} p_\Delta ^2 |\mathcal X(r)|  \cdot (L^\e/r)^2
\lesssim (p_\Delta |\mathcal X(r)|)^2 \cdot L^{-2+(2+\Theta)\e}
= L^{-2+(2+\Theta)\e} \bigl( \bE[|Y|] \bigr)^2\,,
\]
which gives what we needed with a large margin.

\medskip
\paragraph{\underline{The terms with $d_L(x, y)\ge L^\e$}.}
In this case,
\[
p_\Delta(x, y) - p_\Delta^2 \le p(x, y) - p_\Delta^2
= (p(x, y) - p^2) + (p^2-p_\Delta^2)\,.
\]
By the long-distance estimate~\eqref{eq:2point_long} of the two-point function $p(x, y)$,
\[
p(x, y) - p^2 \lesssim L^{-\gamma\e}p^2 \lesssim L^{-\gamma\e} p_\Delta^2\,,
\]
while, by Lemma~\ref{Hessian_conditional},
\[
p^2-p_\Delta^2 \lesssim (\log\Delta)^{1/4} \Delta^{-1/2} p^2
\lesssim L^{-\e/10} p_\Delta^2\,.
\]
Thus, $ p_\Delta(x, y) - p_\Delta^2 \lesssim L^{-c\e}p_\Delta^2$, and the whole double sum
is bounded by $L^{-c\e}p_\Delta^2 \cdot |\mathcal X(r)|^2 =
L^{-c\e}\bigl( \bE[|Y|] \bigr)^2$.

\medskip
This completes the proof of the estimate of $\Var[\, |Y|\, ]$ and hence of
Lemma~\ref{lemma_Cr(alpha)}. \hfill $\Box$

\section{Asymptotic independence}

\begin{lemma}[asymptotic independence]\label{lemma2}
Let $(f_L)$ be a regular Gaussian ensemble, and let
$Z\subset \bS^2(L)$ be an $L^{1-\kappa}$-separated set with sufficiently small positive
$\kappa$.
Then there exist a collection $(\xi(z))_{z\in Z}$ of independent standard
Gaussian random variables and positive constants $c_1, c_2$ so that, for $L\ge L_0$,
\[
\bP \Bigl\{ \max_{z\in Z} | f_L(z) - \xi(z) | > L^{-c_1} \Bigr\}
< e^{-L^{c_2}}\,.
\]
\end{lemma}

\noindent{\em Proof}:
We will follow rather closely the proof of \cite[Theorem~3.1]{NS-IMRN}.
We fix sufficiently large $L$ and introduce the following notation:
\begin{itemize}
\item $\cH_{f_L}$ is a Gaussian Hilbert space generated by $f_L$, i.e.,
the closure of finite linear combinations $\sum c_j f_L(x_j) $ with the scalar product
generated by the covariance.
\item $\cH$ is ``a big Gaussian Hilbert space'' that contains $\cH_{f_L}$ and countably
many mutually orthogonal one-dimensional subspaces that are orthogonal to
$\cH_{f_L}$.
\item $J(z) = f_L(z)$, $z\in Z$, are unit vectors in $\cH_{f_L}$ with
\[
\bigl| \langle J(z), J(z') \rangle_{\cH} \bigr|
= \bigl| \bE[f_L(z) \overline{f_L(z')}] \bigr|
\lesssim L^{-\gamma(1-\kappa)}, \qquad z\ne z'.
\]
\end{itemize}

We claim that {\em there exists a collection of
orthonormal vectors $\bigl\{ \widetilde{J}(z)\bigr\}_{z\in Z}\subset\cH$,
such that
\[
\max_{z\in Z}\, \| \widetilde{J}(z) - J(z)\|_{\cH} \le 2L^{-0.45 \gamma}\,.
\]
}

To prove this claim,
we consider the Hermitian matrix $\Gamma=\Gamma(z, z')_{z, z'\in Z}$ with the elements
\[
\Gamma (z, z') =
\begin{cases}
- \langle J (z), J (z')\rangle_{\cH}, & z' \ne z, \\
L^{-0.9 \gamma}, & z'=z \,.
\end{cases}
\]
For $L>L_0(\gamma, \kappa)$, the matrix $\Gamma$ is positive-definite.
Indeed, by the classical Gershgorin theorem, each eigenvalue of $\Gamma$ lies in one of the intervals
$( \Gamma(z, z)-t(z), \Gamma(z,z)+t(z) )$ with
\[
t(z) = \sum_{z'\in Z\setminus\{z\}}\, |\Gamma(z,z')|,
\]
so we need to check that, for each $z\in Z$, $t(z)<\Gamma(z,z)$, which is easy to see since
\[
t(z) < |Z| \cdot L^{-\gamma(1-\kappa)} \stackrel{|Z|\lesssim L^{2\kappa}}\lesssim
L^{2\kappa} \cdot L^{-\gamma(1-\kappa)} \ll L^{-0.9\gamma}\,,
\]
provided that $\kappa$ is sufficiently small. Since the Hermitian matrix $\Gamma$
is positive-definite, we can find a collection of vectors
$\{I(z)\}_{z\in Z} \subset \cH \ominus \operatorname{span}\{J(z)\}_{z\in Z}$
with the Gram matrix $\Gamma$.
Note that $\| I(z) \|_{\cH} =\sqrt{\Gamma (z,z)} = L^{-0.45 \gamma}$.
Then we let
\[
\widetilde{J}(z) = \frac{J(z)+I(z)}{\| J(z) + I(z)\|_\mathcal H}\,, \qquad z\in Z.
\]
By construction, the system of vectors $\{ \widetilde{J}(z)\}_{z\in Z}$ is orthonormal in $\mathcal H$. Furthermore,
\[
1 \le  \| J(z) + I(z)\|_\mathcal H  \le 1 + L^{-0.45\gamma},
\]
whence,
\[
\| J(z) - \widetilde{J}(z) \|_\mathcal H
\le \bigl( \| J(z) + I(z)\|_\mathcal H - 1 \bigr)\| J(z) \|_{\cH} + \| I(z) \|_{\cH}
\le 2L^{-0.45 \gamma},
\]
proving the claim.

\medskip
It remains to note that, for each $z\in Z$, we have
\[
\cP\bigl\{|\widetilde{J}(z)-J(z)|>t \bigr\} =
e^{-\frac12 (t\|\widetilde{J}(z)-J(z)\|_{\cH}^{-1})^2}
\le e^{-\frac18 t^2 L^{0.9\gamma}},
\]
whence, by the union bound,
\[
\cP\bigl\{ \max_{z\in Z}\, |\widetilde{J}(z)-J(z)|>t \bigr\}
\le |Z|\, e^{-\frac18 t^2 L^{0.9\gamma}}
\lesssim L^{2\kappa} e^{-\frac18 t^2 L^{0.9\gamma}}\,.
\]
Letting, for instance,
$t=L^{-0.4\gamma}$, we complete the proof of Lemma~\ref{lemma2}.
\hfill $\Box$

\section{Two simple lemmas}

In this section we present two simple and standard lemmas which will be employed
later. To keep this work relatively self-contained, we will include their proofs.

\subsection{Anticoncentration of the sums of Bernoulli random variables}

\begin{lemma}\label{lemma:Bernoulli}
Given $0<p_0\le \tfrac12$,
let $(\eta_j)$ be a collection of $N$ independent random variables
on the probability space $\Omega$
such that $\eta_j$ attains the value $1$ with probability $p_j$, $p_0\le p_j\le 1-p_0$, and
the value $0$ with probability $1-p_j$. Let
\[
S_N=\sum_{j=1}^N \eta_j.
\]
Then there exists $\e = \e(p_0)>0$
such that for any measurable function $Q\colon \Omega\to [0, 1]$
with $\displaystyle \int_\Omega Q\, {\rm d}\bP \ge 1-\e $,
and for any $m\in\bR$, we have
\[
\int_{\Omega} |S_N - m|^2 Q\, {\rm d}\bP \ge c(p_0) N.
\]
\end{lemma}

\medskip\noindent{\em Proof of Lemma~\ref{lemma:Bernoulli}}:
During the proof, the value of $c(p_0)$ may vary from
line to line. Take $\la=1/\sqrt{N}$. First, we claim that
$|\bE[\, e^{{\rm i}\la S_N}\,]|\le 1-c(p_0)$.
Indeed,
\[
|\bE[\, e^{{\rm i}\la \eta_j}\,]| = |e^{{\rm i}\la} \cdot p_j + 1\cdot (1-p_j)|
\le 1-c(p_0)\la^2\,,
\]
whence
\[
\bigl| \bE[\, e^{{\rm i}\la S_N}\,] \bigr| \le (1-c(p_0)\la^2)^N \le e^{-c(p_0)N\la^2}
\le 1-c(p_0).
\]
Therefore, for any $m\in\bR$,
\[
\bigl| \bE[\, e^{{\rm i}\la (S_N-m)} - 1 \,]\bigr| \ge c(p_0).
\]
Next, using that $|e^{{\rm i}t}-1|\le |t|$, we proceed as follows:
\begin{align*}
c(p_0) &\le \bigl| \bE[\, (e^{{\rm i}\la (S_N-m)}-1)\cdot (Q+(1-Q))\,]\bigr| \\
&\le \int_{\Omega} \la |S_N-m|\cdot Q\, {\rm d}\bP + \int_{\Omega} 2 (1-Q) {\rm d}\bP
\\
&\le \sqrt{\frac1{N}\,\int_{\Omega} |S_N-m|^2\cdot Q\, {\rm d}\bP} + 2\e.
\end{align*}
Taking $\e\le c(p_0)/4$, we complete the proof.
\hfill $\Box$

\subsection{Large sections lemma}

\begin{lemma}\label{lemma:large_sections}
Let $(\Omega_1 \times \Omega_2, \bP_1 \times \bP_2) $
be a product probability space, and let $0<p\le 1$ and $0<\e \le \tfrac12\, p$.
Let $Q\colon \Omega_1 \times \Omega_2 \to [0, 1]$ be a measurable function with
$\displaystyle \int_{\Omega_1 \times \Omega_2} Q\, {\rm d}\bP \ge 1-\e $,
and let $X\subset \Omega_1$ be an event with
$\bP_1(X)\ge p$. Then
\[
\bP_1 \Bigl\{\omega_1\in X\colon \int_{\Omega_2} Q(\omega_1, \omega_2)\,
{\rm d}\bP_2(\omega_2) \ge 1 - 2\e p^{-1} \Bigr\} \ge \frac{p}2\,.
\]
\end{lemma}

\medskip\noindent{\em Proof of Lemma~\ref{lemma:large_sections}}:
Put
\[
X' = \Bigl\{\omega_1\in X\colon \int_{\Omega_2} Q(\omega_1, \omega_2)\,
{\rm d}\bP_2(\omega_2) \ge 1 - 2\e p^{-1} \Bigr\}\,.
\]
Then,
\begin{align*}
1-\e &\le \int_{\Omega_1 \times \Omega_2} Q\, {\rm d}\bP \\
&= \Bigl( \int_{\Omega_1\setminus X} + \int_{X'} + \int_{X\setminus X'} \Bigr)
\Bigl( \int_{\Omega_2} Q(\omega_1, \omega_2)\,
{\rm d}\bP_2(\omega_2)\Bigr) \, {\rm d}\bP_1(\omega_1) \\
&< 1 - \bP_1(X) + \bP_1(X') + (1-2\e p^{-1})(\bP_1(X)-\bP_1(X')) \\
&\le 1 - 2\e  + 2\e p^{-1}\bP_1(X'),
\end{align*}
which yields the lemma. \hfill $\Box$

\medskip In what follows we will apply this lemma, mostly, with $X=\Omega_1$ and $p=1$.

\section{Variance of the number of loops}

Let $G=G(V, E)$ be a finite graph embedded in the sphere $\bS^2$ with
each vertex having degree four.
We allow $G$ to have multiple edges as well as``circular edges'', which connect a
vertex with itself.
The vertices of the graph are the joints (see
Section~\ref{subsect:saddle-points}), the edges are curves on $\bS^2$ connecting the vertices,
and the faces are the connected components of the open set $\bS^2\setminus (V\cup E)$.

Each vertex $v\in V$ can be replaced by one of two possible ``avoided crossings'' at $v$:
\begin{figure}[h]
\includegraphics[width=0.4\textwidth]{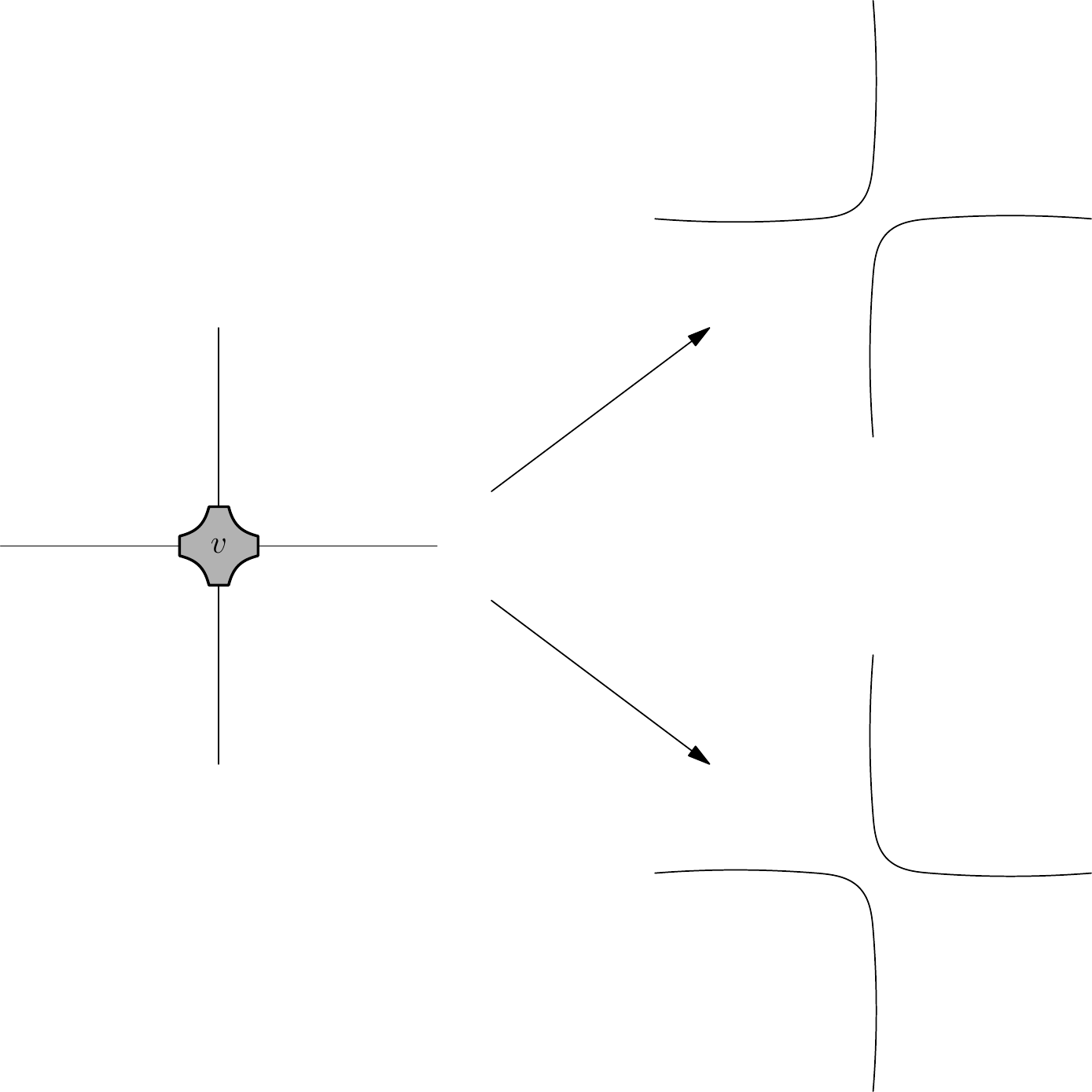}
\caption{The vertex $v$ and its states $\sigma_v$}
\label{zfigureI}
\end{figure}

We call the choice of the avoided crossing at $v$
{\em the state of the vertex $v$} and denote it
by $\sigma_v$.
When the states are assigned to all vertices in $V$, the collection of states $\sigma_V = \{\sigma_v\colon v\in V\}$ turns the graphs $G$ into a collection of loops $\Gamma=\Gamma(\sigma_V)$.

We will deal with a random loop model when the states $\sigma_v$ are independent
random variables taking their values with probabilities $p(v)$ and $1-p(v)$.
By $(\Omega, \bP)$ we denote the probability space
on which the random states are defined. Then $N(\Gamma)=N(\Gamma(\sigma_V))$
is a random variable on $(\Omega, \bP)$.

Given $0<p_0\le\frac12$, we put $V(p_0)=\bigl\{v\in V\colon p_0 \le p(v) \le 1-p_0 \bigr\}$
and denote by $|V(p_0)|$ the cardinality of the set of vertices $V(p_0)$.

\begin{lemma}\label{lemma:variance_clusters}
For any $0<p_0\le\frac12$, there exist positive $c(p_0)$, $C(p_0)$, and $\e=\e(p_0)$
such that for any function $Q(\sigma_V)$ defined on the set of all possible states
and taking the values in the interval $[0, 1]$ with
$\displaystyle \int_\Omega Q(\sigma_V)\, {\rm d}\bP  \ge 1-\e$, and for any $m\in\bR$,
\[
\int_{\Omega} ( N(\Gamma(\sigma_V)) - m)^2\, Q(\sigma_V)\,
{\rm d}\bP \ge c(p_0) |V(p_0)|\,,
\]
provided that $|V(p_0)|\ge C(p_0)$.
\end{lemma}

\subsection{Beginning the proof of Lemma~\ref{lemma:variance_clusters}}
We fix a function $Q$ as above.
In several steps we will reduce the statement of the lemma to the anti-concentration bound for the sum of independent Bernoulli random variables provided by Lemma~\ref{lemma:Bernoulli}.
In each of these steps we will be using the following decoupling argument.

\subsubsection{Decoupling}\label{subsubsect:decoupling}

Suppose that the vertices are split into two disjoint parts:
$ V = V' \sqcup V'' $ and decompose correspondingly
$\sigma_V = (\sigma_{V'}, \sigma_{V''})$.
A collection of states $\sigma_{V'}$ assigned to the vertices from $V'$ generates
\begin{itemize}
\item a collection $\Gamma'=\Gamma (\sigma_{V'})$ of disjoint loops,
\item and a graph
$G(\sigma_{V'})$
with vertices at $V''$, all of them having degree $4$.
\end{itemize}
\begin{figure}[h]
\includegraphics[width=0.6\textwidth]{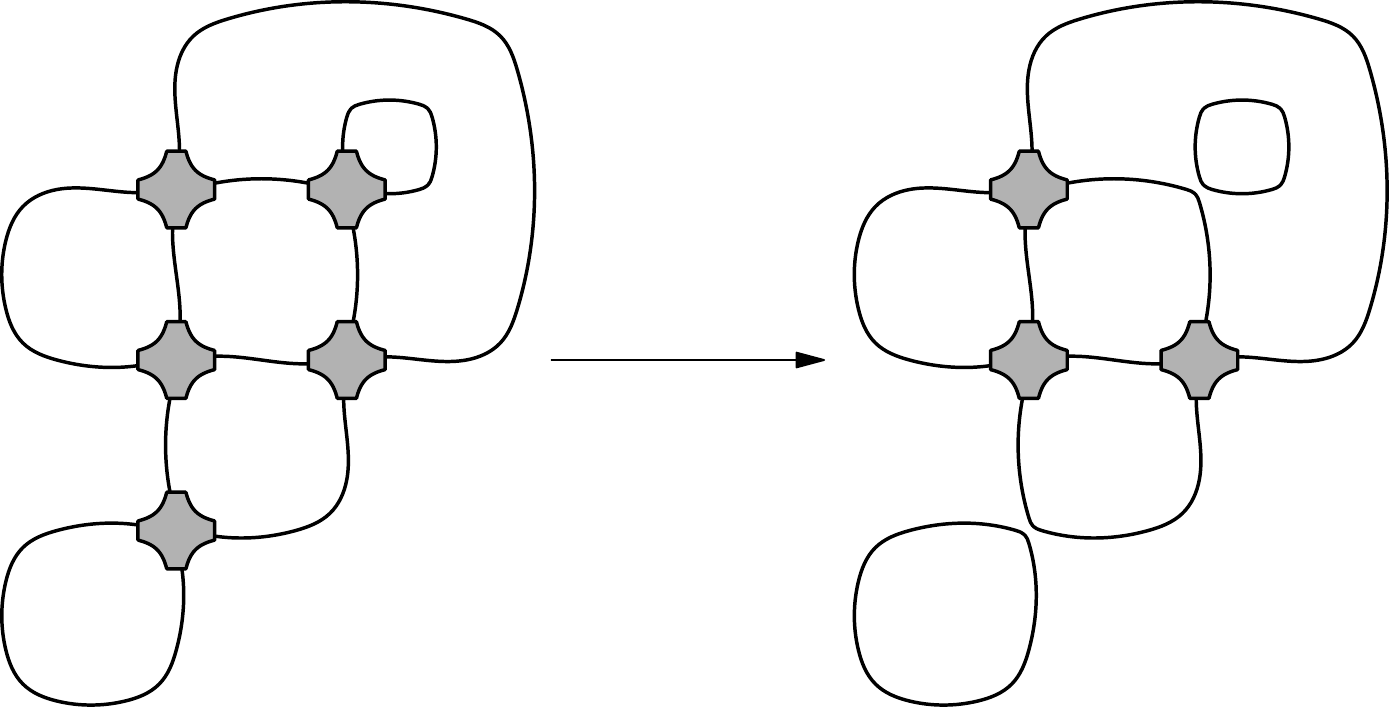}
\caption{On the left: the graph $G$. On the right: the loops $\Gamma'$ and the graph $G(\sigma_{V'})$.}
\label{zfigureII}
\end{figure}
By $N(\Gamma')$ we denote the number of loops in the collection $\Gamma'$.
The collection of states $\sigma_{V''}$ turns the graph $G(\sigma_{V'})$
into a collection of loops $\Gamma'' = \Gamma (\sigma_{V'}, \sigma_{V''})$.
Then,
$\Gamma = \Gamma' \sqcup \Gamma''$, and
$N(\Gamma) = N(\Gamma') +
N(\Gamma'')$.

Let $\phi (\sigma_V) $ be any non-negative bounded measurable function.
Since the random variables $ \sigma_{V'}  $ and
$ \sigma_{V''} $ are independent, we have
\begin{multline*}
\int_\Omega \phi (\sigma_V(\omega))\, {\rm d}\bP (\omega)
= \iint_{\Omega\times\Omega}
\phi(\sigma_{V'}(\omega'), \sigma_{V''}(\omega''))\,
{\rm d}\bP(\omega') {\rm d}\bP(\omega'')\,.
\\
= \int_{\Omega} \Bigl[
\int_{\Omega}
\phi(\sigma_{V'}(\omega'), \sigma_{V''}(\omega''))\,
{\rm d}\bP(\omega'') \Bigr]\, {\rm d}\bP(\omega')\,.
\end{multline*}
Hence, for any event $X'\subset \Omega$,
\[
\int_\Omega \phi (\sigma_V)\, {\rm d}\bP
\ge \bP (X') \cdot \inf_{\sigma_{V'}\in\Sigma'}\,
\int_{\Omega}
\phi(\sigma_{V'}, \sigma_{V''}(\omega''))\,
{\rm d}\bP(\omega'')\,,
\]
where $\Sigma' = \{\sigma_{V'}(\omega')\colon \omega'\in X'\}$.
Letting $\phi = (N(\Gamma)-m)^2\, Q$ and  taking into account
that $ N(\Gamma) - m = N(\Gamma'') - (m - N(\Gamma'))$, we get
\begin{multline*}
\int_\Omega (N(\Gamma(\sigma_V))-m)^2 Q(\sigma_V)\, {\rm d}\bP \\
\, \ge \, \bP(X')\, \cdot \inf_{\sigma_{V'}\in\Sigma'}\, \inf_{\ell \in\bR}\,
\int_{\Omega} (N(\Gamma''(\sigma_{V''}(\omega''))) - \ell)^2\,
Q(\sigma_{V'}, \sigma_{V''}(\omega''))\, {\rm d}\bP(\omega'')\,.
\end{multline*}
The choice of the event $X'\subset\Omega$, or what is the same, of the set of
states $\Sigma'$,
is in our hands. Choosing it we need to keep the value of the integral
\[
\inf_{\sigma_{V'}\in\Sigma'}\, \int_{\Omega} Q(\sigma_{V'}, \sigma_{V''}(\omega''))\,
{\rm d}\bP (\omega'')
\]
close to $1$, while $ \bP (X') $ should stay bounded away from zero.
This will be done with the help of Lemma~\ref{lemma:large_sections}.
After that it will suffice to prove Lemma~\ref{lemma:variance_clusters}
for the graph $G(\sigma_{V'})$ with vertices at $V''$.

In what follows, we will apply this decoupling argument several times.
To simplify the notation, after each step we treat the function $Q$
as depending only on the states $\sigma_{V''}$ of the remaining set of vertices $V''$,
ignoring its dependence on the fixed states $ \sigma_{V'}\in\Sigma'$.

\subsection{Discarding the vertices $v$ with $p(v)<p_0$ or $p(v)>1-p_0$}

As above, we let $V(p_0)=\bigl\{v\in V\colon p_0 \le p(v) \le 1-p_0 \bigr\}$.
We put $V' = V\setminus V(p_0)$, $V'' = V(p_0)$, and consider
\[
X' = \Bigl\{\omega'\in\Omega\colon
\int_\Omega Q(\sigma_{V'}(\omega'),
\sigma_{V''}(\omega'')) \,{\rm d}\bP(\omega'') \ge 1 - 2\e \Bigr\}\,.
\]
Then, by Lemma~\ref{lemma:large_sections} (applied with $p=1$),
$ \bP(X') \ge \tfrac12 $.

So, from now on, we assume that for each vertex $v$ of the graph $G$, we have
$p_0\le p(v) \le 1-p_0$, while
$\displaystyle \int_\Omega Q(\sigma_V)\, {\rm d}\bP \ge 1-2\e$.

\subsection{Many faces have at most $4$ vertices on the boundary}

Denote by $F$ the set of the faces of the graph $G$.
Given a face $\mathfrak f\in F$, we denote by $\mathfrak v(\mathfrak f)$ the set of vertices
in $V$ that lie on the boundary $\partial\mathfrak f$. The Euler formula gives us
\[
|V| - |E| + |F| = 1 + \nu,
\]
where $\nu $ is the number of connected components of the graph $G$.
Since the degree of each vertex in $G$ equals $4$, we have
$|E| = 2 |V| $, whence,
$ |F| \ge |V| + 2 $.
Furthermore, since each vertex in $G$ lies on the boundary of at most four
different faces, we have
\[
|F| > |V| \ge \frac14\, \sum_{\mathfrak f\in F} |\mathfrak v (\mathfrak f)|
\ge \frac14\,
\sum_{\substack{\mathfrak f\in F \\ |\mathfrak v( \mathfrak f)|\ge 5}}
|\mathfrak v(\mathfrak f)|
\ge \frac54\, \bigl| \{\mathfrak f\in F\colon |\mathfrak v( \mathfrak f)|\ge 5\} \bigr|.
\]
We let $F^* = \{\mathfrak f\in F\colon |\mathfrak v( \mathfrak f)|\le 4\}$ be the set of all faces having at most $4$ vertices on the boundary and conclude that
$ \bigl| F^* \bigr| > \tfrac15 |F| > \tfrac15 |V| $.

\subsection{Choosing a maximal collection of separated faces}
We call faces $\mathfrak f, \mathfrak f'\in F^*$ {\em separated\ } if they do not
have common vertices on their boundaries:
$\mathfrak v( \mathfrak f) \cap \mathfrak v( \mathfrak f')=\emptyset$.
We fix a maximal collection $\mathfrak F\subset F^*$ of separated faces.
Since for any face $\mathfrak f\in \mathfrak F$ there are at most $12$ other
faces $\mathfrak f'\in F^*$ with $\mathfrak v(\mathfrak f') \cap \mathfrak v(\mathfrak f)
\ne \emptyset$, we conclude from the maximality of $\mathfrak F$ that it also contains sufficiently many faces:
\[
| \mathfrak F | \ge \frac1{13}\, \bigl| F^* \bigr| > \frac1{65} |V|.
\]

\subsection{Marking vertices and cycles}
Next,  we choose a simple cycle $\mathfrak c$ on the boundary of each face
$\mathfrak f\in\mathfrak F$, and mark a vertex
$v_\mathfrak c$ on that cycle according to the following rule:
We take a vertex $v'$ in $\mathfrak v(\mathfrak f)$ and, starting
at $v'$, walk along $\partial\mathfrak f$ turning left at each vertex
so that the face $\mathfrak f$ always remains on the left-hand side.
We stop when we return for the first time to the vertex that we already have passed,
and mark that vertex and the corresponding cycle.

We will be using the following property of the marked cycles: {\em if the cycle exits a vertex along an edge $e$, then it returns to this vertex along an edge, which is one of two edges adjacent to $e$} (this follows from the fact that the same face cannot lie on both sides of some edge). This property yields that there exist states of vertices on $\mathfrak c$ which turn $\mathfrak c$ into a separate loop.

We denote by $V_{\tt M}\subset V$  the set of all marked vertices, and
by $V_{\tt UM}=V\setminus V_{\tt M}$ the set of all unmarked vertices.

\subsection{Good marked cycles}
For any particular assignment of the states of the unmarked vertices,
the marked cycle $\mathfrak c$ is called {\em good} if the edges of $\mathfrak c$
merge into an
edge of the graph $G(\sigma(V_{\tt UM}))$ that connects the vertex $v_{\mathfrak c}$
with itself.
\begin{figure}[h]
\includegraphics[width=0.5\textwidth]{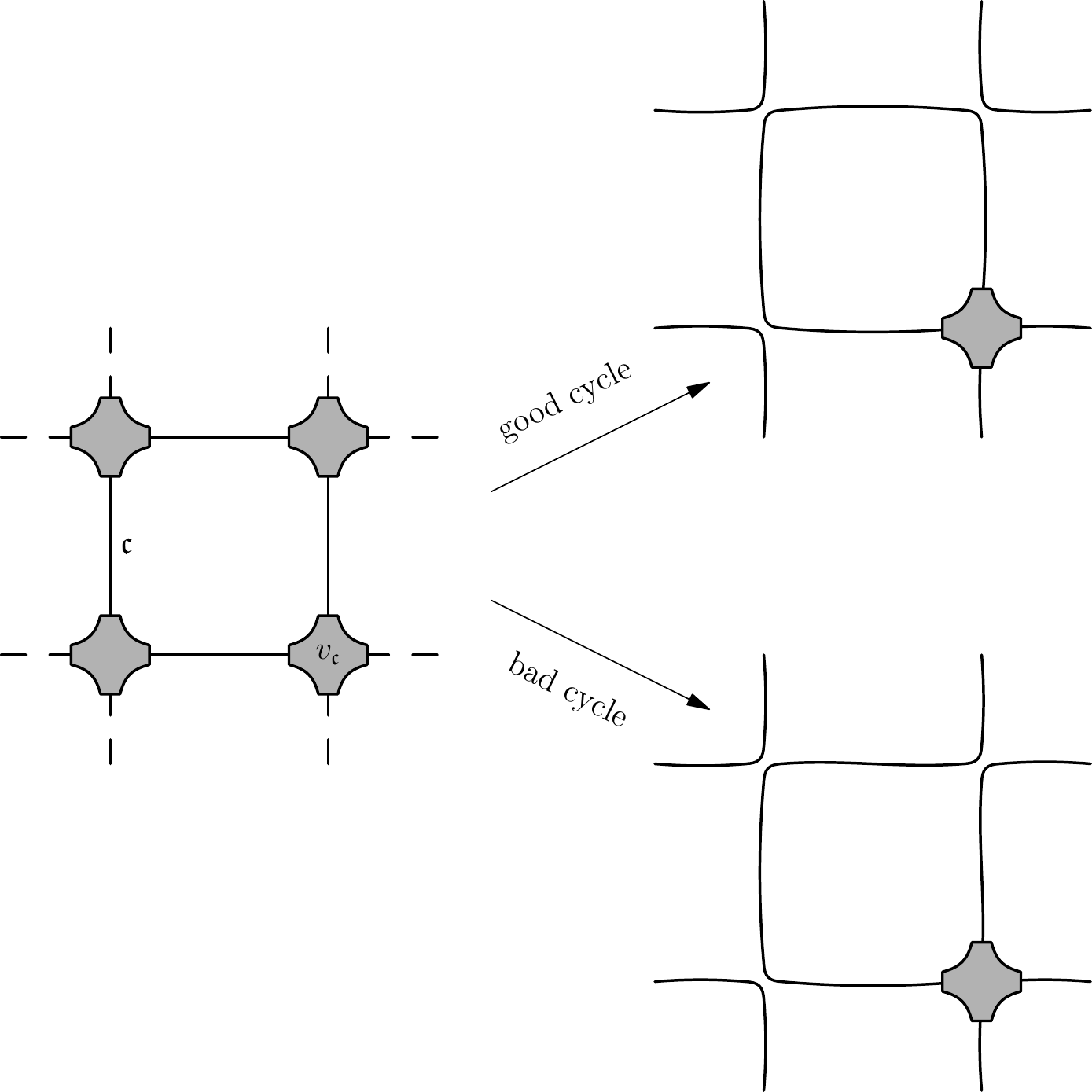}
\caption{Good and bad cycles}
\label{zfigureVI}
\end{figure}
Note that this event
depends only on the states of
at most $3$ unmarked vertices lying on $\mathfrak c\setminus\{v_\mathfrak c\}$.
Therefore, for any marked cycle $\mathfrak c$, we have
$ \bP \bigl[\, \mathfrak c \ {\rm is\ good\ } \bigr] \ge p_0^3 $.
Hence, denoting by $N_{\tt G}$
the number of good cycles, we obtain
\[
\bE \bigl[ N_{\tt G} \bigr] \ge p_0^3\, |\{{\rm marked\ cycles}\}|
= p_0^3\, |\mathfrak F| \ge \frac1{65}\, p_0^3\, |V|.
\]
Using first the Chebyshev inequality and then the independence of the random
states, we get
\begin{align*}
\bP \bigl\{ N_{\tt G} \le \tfrac12\, \bE [N_{\tt G}] \bigr\}
&\le \frac{4\operatorname{Var} [N_{\tt G}]}{\bigl( \bE [N_{\tt G}] \bigr)^2} \\
&\le 130^2 p_0^{-6}\, \frac{|\{{\rm marked\ cycles}\}|}{|V|^2} \\
&\le 130^2 p_0^{-6}\, |V|^{-1}.
\end{align*}
To simplify the notation, we let $ d(p_0)=\tfrac1{130}\, p_0^3 $ and $D(p_0)=d(p_0)^{-2}$.
Then, letting $ X = \bigl\{\omega'\colon N_{\tt G} > d(p_0)|V| \bigr\}$,
we see that $ \bP (X) \ge p \stackrel{\rm def}= 1- D(p_0)|V|^{-1} $.
Put
\[
X' = \Bigl\{\omega'\in X\colon
\int_\Omega Q(\sigma_{V_{\tt UM}}(\omega'), \sigma_{V_{\tt M}}(\omega''))\,
{\rm d}\bP(\omega'') \ge 1 - 4\e p^{-1} \Bigr\}
\]
(recall that at this moment $\displaystyle \int_\Omega Q\, {\rm d}\bP \ge 1-2\e$).
Then, by Lemma~\ref{lemma:large_sections},
$\bP (X') \ge \tfrac12 p$, which is $ \ge \tfrac14$, provided that $|V|> 2D(p_0)$.
From now on, we fix the states of the unmarked vertices corresponding to the event $X'$
and consider the remaining graph with vertices in $V_{\tt M}$.
At this step, $\displaystyle \int_\Omega Q\, {\rm d}\bP \ge 1-4\e p^{-1}$.

\subsection{Discarding bad cycles}
All marked cycles $\mathfrak c$ are split into two classes:
bad cycles and good cycles. Correspondingly, we decompose the set of all marked
vertices $V_{\tt M}$ into the disjoint union $V_{\tt M} = V_{\tt M.B.} \sqcup V_{\tt M.G.}$, and consider
\[
X' = \Bigl\{\omega'\in\Omega\colon \int_\Omega Q(\sigma_{V_{\tt M.B.}}(\omega'),
\sigma_{V_{\tt M.G.}}(\omega''))\,
{\rm d}\bP(\omega'') \ge 1  - 8\e p^{-1} \Bigr\}
\]
By Lemma~\ref{lemma:large_sections} (applied with $X=\Omega$),
$\bP(X')\ge \tfrac12$. We fix the states of marked bad vertices
corresponding to the event $X'$.

\subsection{Completing the proof of Lemma~\ref{lemma:variance_clusters}}
We are left with the graph $G$ with vertices at $V_{\tt M.G.}$.
For each vertex $v\in V_{\tt M.G.}$,
there is ``a circular edge'' $e_v$ with the endpoints
at $v$, which came from the corresponding cycle $\mathfrak c$.
Each state $\sigma_v$ of the vertex $v$
either creates from this circular edge a separate loop, or merges it with other egdes:
\begin{figure}[h]
\includegraphics[width=0.5\textwidth]{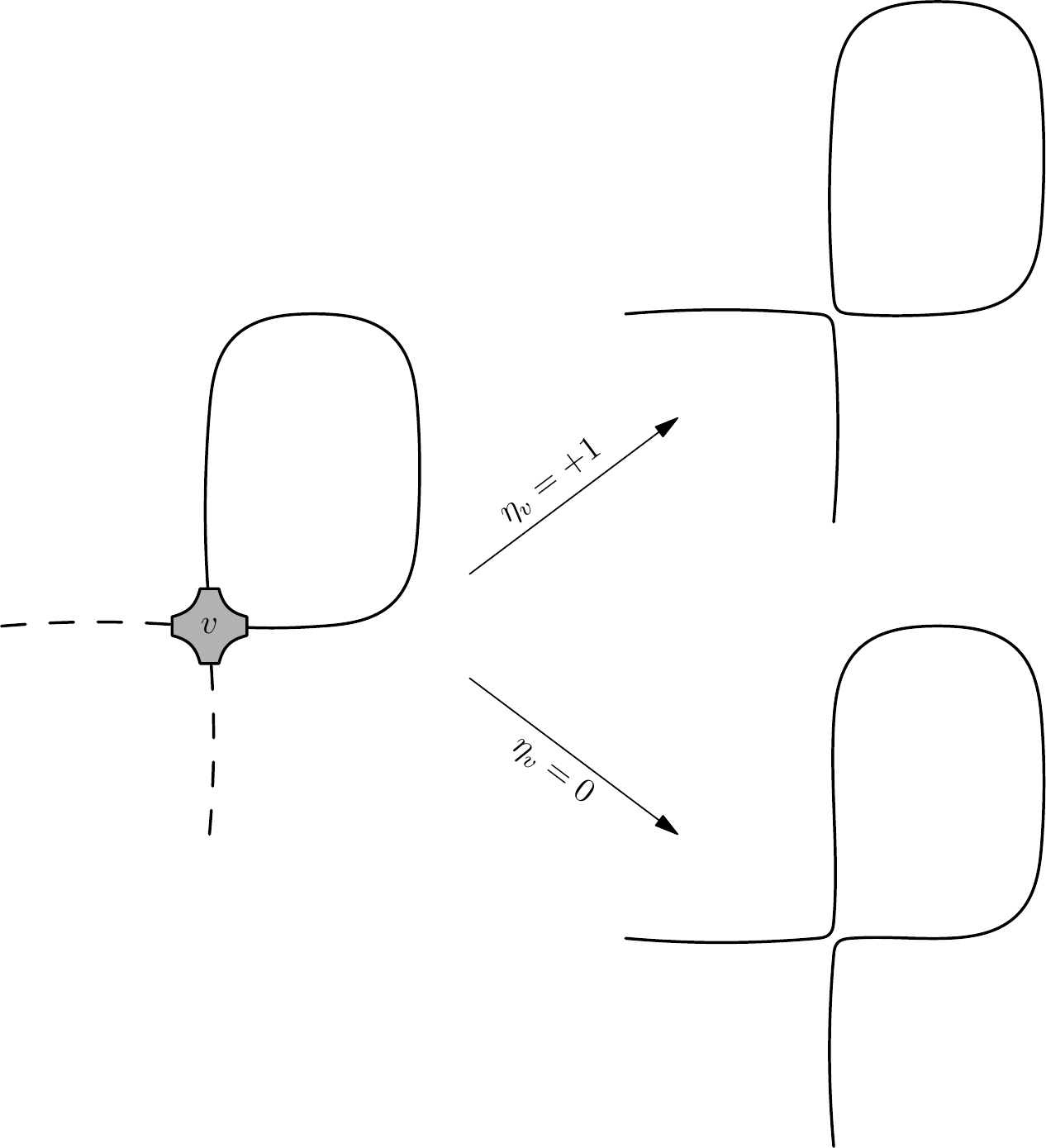}
\caption{Separation vs merging}
\label{zfigureVII}
\end{figure}

We fix a collection of states
$\sigma^*_{V_{\tt M.G.}}$ of good marked vertices such that
that none of the corresponding good cycles turns into a separate loop, and denote by
$\Gamma( \sigma^*_{V_{\tt M.G.}} )$ the loop ensemble
obtained from the graph $G$ after the assignment of the states
$\sigma^*_{V_{\tt M.G.}}$.
Introduce a collection of independent Bernoulli random variables
$\bigl\{ \eta_v \bigr\}_{v\in V_{\tt M.G.}}$, letting
$\eta_v=0$ if $\sigma_v = \sigma^*_{V_{\tt M.G.}}(v)$,
and $\eta_v=1$ otherwise. Then
\[
N(\Gamma(\sigma_{V_{\tt M.G.}})) = N(\Gamma(\sigma^*_{V_{\tt M.G.}}))
+ \sum_{v\in V_{\tt M.G.}} \eta_v\,.
\]
Applying Lemma~\ref{lemma:Bernoulli}, we get the uniform in $m$ lower bound
\[
\int_{\Omega}
\bigl( N(\Gamma(\sigma_{V_{\tt M. G.}})) - m \bigr)^2\, Q (\sigma_{V_{\tt M.G.}})\, {\rm d}\bP
\ge c(p_0)\, |V_{\tt M.G.}|\,.
\]
To finish off the proof of Lemma~\ref{lemma:variance_clusters}, it remains to recall that
good marked vertices are in the one-to-one correspondence with good cycles, that is,
$ | V_{\tt M.G.} | = N_{\tt G} $, and that the
states of unmarked vertices were fixed so that
$ N_{\tt G} \ge d(p_0) |V| $. \mbox{} \hfill $\Box$

\section{Tying loose ends together: proof of the theorem}

\subsection{Perturbing $f_L$}
We choose a sufficiently small $\eps>0$, take $\alpha'=L^{-2+\eps}$, $\alpha=L^{-2+2\eps}$.
Then we take the function $f_L$ and its independent copy $g_L$, and put
\[
\widetilde{f}_L = \sqrt{1-\alpha'^2} f_L + \alpha' g_L.
\]
This is a random Gaussian function equidistributed with $f_L$.
We will show that
\[
\inf_{m\in\bR} \bE\bigl[ (N(\widetilde{f}_L) -m)^2 \bigr] \gtrsim L^{c\e},
\]
which immediately yields the lower bound for $\Var[N(f_L)] $ we are after.
Note that
\[
\bE\bigl[ (N(\widetilde{f}_L) -m)^2 \bigr]
= \bE^{f_L}\, \bE^{g_L}
\bigl[ N(\sqrt{1-\alpha'^2}f_L +\alpha' g_L) - m )^2\bigr].
\]
That is, it suffices to show that with probability at least $\tfrac12$ in $f_L$ we have
\[
\bE^{g_L} \bigl[ (N(\widetilde{f}_L) -m)^2 \bigr] \gtrsim L^{c\e}.
\]

\subsection{Freezing $f_L$}
We will prove a somewhat stronger statement that this inequality holds if the function $f_L$
satisfies the following conditions:
\begin{itemize}
\item $f_L\in C^3(A, \Delta, \alpha, \beta)$ (introduced in~\ref{subsect:sets_Cr_class_C3}) with $A=\log L$, $\Delta=L^{3\e}$
and with $\beta$ chosen so that
$\beta^2 L^{7\e} = \alpha$ (i.e., $\beta=L^{-1-\frac52\e}$);
\item the set $\Cr(\alpha)$ is $L^{1-C\eps}$-separated;
\item $|\Cr(\alpha')|\ge L^{c\eps}$.
\end{itemize}
By Lemma~\ref{lemma-Hessian-lowerbound1} and Lemma~\ref{lemma_Cr(alpha)}, these
three conditions hold w.h.p. in $f_L$. From now on, we fix $f_L$ so
that these conditions hold, and omit the index $g_L$ meaning $\bP=\bP^{g_L}$,
$\bE=\bE^{g_L}$, etc.

\subsection{Recalling a little Morse caricature}
The rest will essentially follow from our little Morse caricature summarized in Lemma~\ref{LemmaM} combined with Lemma~\ref{lemma:variance_clusters} on the fluctuations
in the number of random loops. In order to apply Lemma~\ref{LemmaM}, first, we observe
that the relations
\[
A\alpha' \ll \alpha, \quad A\Delta^2\beta^2 \ll \alpha \ll (A\Delta)^{-2}\beta,
\quad A^2 \Delta^3\alpha \ll 1
\]
required in Lemma~\ref{LemmaM} readily follow from our choice of the parameters
$\alpha$, $\alpha'$, $\beta$, $A$ and $\Delta$ made few lines above.
Lemma~\ref{LemmaM} also needs the lower bound
\[
\min_{\Cr(\alpha)} |\widetilde{f}_L| \gtrsim A \Delta^2\alpha^2,
\]
which holds w.h.p. with a large margin since $A \Delta^2\alpha^2 = L^{-4+10\e}\, \log L $ while,
as we will momentarily see, w.h.p. in $g_L$, we have
\begin{equation}\label{eq:*a}
\min_{\Cr(\alpha)} |\widetilde{f}_L| > \alpha' L^{-c_1}\,,
\end{equation}
where $c_1\le 1$ is a constant from Lemma~\ref{lemma2} (recall that $\alpha' = L^{-2+\e}$).
Indeed, since $g_L(p)$ is a standard Gaussian random variable, the probability that
\[
| \sqrt{1-\alpha'^2} f_L(p) + \alpha' g_L (p) | \le \alpha' L^{-c_1}
\]
at a given point $p$ is $\lesssim L^{-c_1}$. By the union bound, the probability
that this happens somewhere on $\Cr(\alpha)$ is
\[
\lesssim L^{-c_1} |\Cr(\alpha)| \lesssim L^{-c_1+2C\e}
\ll L^{-c_1/2},
\]
provided that $\eps$ is sufficiently small.

Thus, Lemma~\ref{LemmaM} applied to the
functions $f_L$ and $\widetilde{f}_L$ yields that
\[
N(\widetilde{f}_L) = N_{\rm I}(\widetilde{f}_L) + N_{\rm II}(\widetilde{f}_L)
+ N_{\rm III}(\widetilde{f}_L)
\]
on the major part of the probability space where $g_L\in C^3(\log L)$ and where estimate~\eqref{eq:*a} holds. The first term on the RHS, $N_{\rm I}(\widetilde{f}_L)$, comes from the stable connected components of $Z(f_L)$. Hence, on the large part of the probability space, the fluctuations
in $N(\widetilde{f}_L)$ come only from the blinking circles
$ N_{\rm II}(\widetilde{f}_L) $ and from the Bogomolny-Schmit loops
$ N_{\rm III}(\widetilde{f}_L) $, and  after we have fixed the function $f_L$,
both these quantities depend only on the configuration of (random) signs
of $\widetilde{f}_L(p)$, $p\in\Cr(\alpha)$.

\subsection{Fight for independence}
To make these random signs independent,
using Lemma~\ref{lemma2}, we choose
a collection of {\em independent}
standard Gaussian random variables $\xi(p)$, $p\in\Cr(\alpha)$, so that
\begin{equation}\label{eq:**}
\max_{p\in\Cr(\alpha)} |g_L(p)-\xi(p)|\le L^{-c_1}\,.
\end{equation}
Denote by $\Omega'$ the event that $\| g_L \|_{C^3} \le \log L$ and
both estimates~\eqref{eq:*a} and~\eqref{eq:**}  hold. Then
\[
\tau(L)\stackrel{\rm def}=\bP(\Omega\setminus\Omega')=o(1),
\qquad L\to\infty,
\]
while on $\Omega'$ we have
\[
\sgn\bigl(\widetilde{f}_L (p)\bigr) = s(p) \stackrel{\rm def}=
\sgn\bigl( \sqrt{1-\alpha'^2}f_L(p) + \alpha' \xi (p) \bigr),
\qquad p\in\Cr(\alpha).
\]
Note that the random signs $s(p)$ are independent and that there exists
$p_0\in (0, \tfrac12]$ such that on $\Cr(\alpha')$ each of the two possible values
of $s(p)$ is attained with probability at least $p_0$.

For any subset $Z\subset \Cr(\alpha)$, we let
$s_{Z} = \bigl( s(p) \bigr)_{p\in Z}$. As before, we use the notation
\begin{align*}
\operatorname{Cr}_{\tt S}(\alpha) &=
\{p\in\Cr(\alpha)\colon p\  {\rm is\ a\ saddle\ point\ of\ } f_L \}, \\
\operatorname{Cr}_{\tt E}(\alpha) &=
\{p\in\Cr(\alpha)\colon p\  {\rm is\ a\ local\ extremum\ of\ } f_L \}.
\end{align*}
Since, on $\Omega'$,
$ N_{\rm II}(\widetilde{f}_L) $ depends only on
$s_{\operatorname{Cr}_{\tt E}(\alpha)}$ and $ N_{\rm III}(\widetilde{f}_L) $
on $s_{\operatorname{Cr}_{\tt S}(\alpha)}$,
there exist functions $\widetilde{N}_{\rm II}(s_{\operatorname{Cr}_{\tt E}(\alpha)})$ and
$ \widetilde{N}_{\rm III}(s_{\operatorname{Cr}_{\tt S}(\alpha)})$ such that, on $\Omega'$,
we have
$ N_{\rm II}(\widetilde{f}_L) =
\widetilde{N}_{\rm II}(s_{\operatorname{Cr}_{\tt E}(\alpha)})$
and
$ N_{\rm III}(\widetilde{f}_L)
= \widetilde{N}_{\rm III}(s_{\operatorname{Cr}_{\tt S}(\alpha)})$.
The function $N_{\rm I}(\widetilde{f}_L)$ stays constant on $\Omega'$, by
$\widetilde{N}_{\rm I}$ we denote the value of that constant.
Denoting by $\chi_{\Omega'}$ the indicator-function of the event $\Omega'$,
we get
\begin{align*}
\bE\bigl[ (N(\widetilde{f}_L) - m)^2 \bigr] &\ge
\int_\Omega [ N_{\rm II}(\widetilde{f}_L) + N_{\rm III}(\widetilde{f}_L) -
(m - N_{\rm I}(\widetilde{f}_L))]^2 \chi_{\Omega'}\, {\rm d}\bP \\
&= \int_\Omega [\widetilde{N}_{\rm II} + \widetilde{N}_{\rm III} -
(m - \widetilde{N}_{\rm I})]^2\, \bE[\chi_{\Omega'}\big| s_{\Cr(\alpha)}]\,
{\rm d}\bP\,.
\end{align*}
The conditional expectation $ \bE[\chi_{\Omega'}\big| s_{\Cr(\alpha)}] $
can be written as $Q(s_{\Cr(\alpha)})$, where $Q$ is a function
on a finite set $S_{\Cr(\alpha)}$ of all possible collections of signs $s_{\Cr(\alpha)}$.
Thus,
\[
\bE\bigl[ (N(\widetilde{f}_L) - m)^2 \bigr] \ge
\int_\Omega [\widetilde{N}_{\rm II} + \widetilde{N}_{\rm III} -
(m - \widetilde{N}_{\rm I})]^2\, Q\, {\rm d}\bP\,.
\]
Note that $\bE [Q] = \bP(\Omega') = 1-\tau(L)$.

\subsection{Fluctuations generated by the Bogomolny-Schmit loops}
First, we consider the case when $|\Cr_{\tt S}(\alpha')|
\ge \tfrac12 |\Cr(\alpha')|$ and look at the fluctuations in
the number of the Bogomolny-Schmit loops. We use a decoupling
argument similar to the one introduced in Section~\ref{subsubsect:decoupling}.
We decompose $s_{\Cr(\alpha)} = (s_{\Cr_{\tt E}(\alpha)}, s_{\Cr_{\tt S}(\alpha)})$
and let
\[
X' = \Bigl\{\omega_1\in\Omega\colon
\int_\Omega Q(s_{\Cr_{\tt E}(\alpha)}(\omega_1), s_{\Cr_{\tt S}(\alpha)}(\omega_2))\,
{\rm d}\bP(\omega_2) \ge 1 - 2\tau(L) \Bigr\}\,.
\]
Then, by Lemma~\ref{lemma:large_sections}, $\bP(X')\ge \tfrac12$,
and therefore,
\begin{multline*}
\bE\bigl[ (N(\widetilde{f}_L) - m)^2 \bigr] \\
\quad \ge \tfrac12\, \underset{\omega_1\in X'}{\operatorname{ess\, inf}}\,
\int_\Omega
\bigl[ \widetilde{N}_{\rm III}(s_{\Cr_{\tt S}}(\omega_2)) -
(m - \widetilde{N}_{\rm I} - \widetilde{N}_{\rm II}(s_{\Cr_{\tt E}}(\omega_1)))\bigr]^2\,
Q(s_{\Cr_{\tt E}}(\omega_1), s_{\Cr_{\tt S}}(\omega_2))\, {\rm d}\bP (\omega_2)\,.
\end{multline*}

We fix $\omega_1\in X'$ and the corresponding signs $s_{\Cr_{\tt E}(\alpha)}(\omega_1)$ and
consider the graph $G(V, E)$ introduced in~\ref{subsubsect:BS-loops}. The vertices of this graph are the joints $J(p, \delta)$ with
$p\in\Cr_{\tt S}(\alpha)$ and $\delta = c(A\Delta)^{-1}$ with sufficiently small
positive constant $c$.
The edges are connected components of the set
\[
Z(f_L)\setminus\bigcup_{p\in\Cr_{\tt S}(\alpha)} J(p, \delta)
\]
that touch the boundaries $\partial J(p, \delta)$. The random states $\sigma_v$
are defined by the signs $s(p)$, and, by construction, are independent.
Furthermore, for the vertices $v$ corresponding to the set $\Cr_{\tt S}(\alpha')$, the probabilities of the states $\sigma_v$ lie in the range $[p_0, 1-p_0]$.
Then Lemma~\ref{lemma:variance_clusters} yields that
\[
\inf_{m\in \bR}\,
\bE \bigl[ (N(\widetilde{f}_L)-m)^2 \bigr]
\gtrsim |\Cr_{\tt S}(\alpha')| \ge \tfrac12\,L^{c\e}.
\]
This finishes the proof of our theorem in the case when
$|\Cr_{\tt S}(\alpha')| \ge \tfrac12 |\Cr(\alpha')|$.

\subsection{Fluctuations generated by blinking circles}
It remains to consider the case when at least half of the critical points in $\Cr(\alpha')$
are local extrema. In this case, we use the decomposition
$\Cr(\alpha) = (\Cr(\alpha)\setminus \Cr_{\tt E}(\alpha')) \sqcup \Cr_{\tt E}(\alpha') $
and once again combine Lemma~\ref{lemma:large_sections} with Lemma~\ref{lemma:Bernoulli}.
Let
\[
X' = \Bigl\{\omega_1\in\Omega\colon
\int_\Omega Q(s_{\Cr(\alpha)\setminus \Cr_{\tt E}(\alpha')}(\omega_1),
s_{\Cr_{\tt E}(\alpha')}(\omega_2))\, {\rm d}\bP(\omega_2) \ge 1 - 2\tau(L) \Bigr\}\,.
\]
Then, by Lemma~\ref{lemma:large_sections}, $\bP(X')\ge\tfrac12$.
We fix $\omega_1\in X'$ and the corresponding
$s_{\Cr(\alpha)\setminus \Cr_{\tt E}(\alpha')}(\omega_1)$.
The value of $\widetilde{N}_{\rm II}$ is the number
of $ p\in\Cr_{\tt E}(\alpha') $ such that the sign $s(p)$
is opposite to the sign of the
eigenvalues of $ \cH_{f_L}(p) $.
To each $p\in \Cr_{\tt E}(\alpha')$ we associate a Bernoulli random variable
\[
\eta_p =
\begin{cases}
1 & s(p) \
{\rm is\ opposite\ to\ the\ sign\ of\ the\ eigenvalues\ of\ } \cH_{f_L}(p), \\
0 & {\rm otherwise}.
\end{cases}
\]
Since the signs $s(p)$ are independent, the variables $\eta_p$ are independent as well.
Recall that
everywhere on $\Cr(\alpha')$ each of two possible values of $s(p)$
is attained with probability $\ge p_0$,
and note that
\[
\widetilde{N}_{\rm II} = \sum_{p\in \Cr_{\tt E}(\alpha')} \eta_p.
\]
Then, Lemma~\ref{lemma:Bernoulli} does the job. This
finishes off the proof of our theorem in the second case when
$|\Cr_{\tt E}(\alpha')| \ge \tfrac12\, |\Cr (\alpha')| $. \hfill $\Box$

\section{The case of spherical harmonics}

As we have already mentioned in the introduction, our theorem does not straightforwardly
apply to the ensemble of Gaussian spherical harmonics.
Here, we will outline minor modifications needed in this case.

\subsection{}
The spherical harmonic $f_n$ is an even (when its degree $n$ is even) or odd (when $n$
is odd) function.
Hence, its zero set $Z(f_n)$ is symmetric with respect to the origin.
So, for the advanced readers, we are just working
on the projective space $\bR\bP^2$ instead of the sphere. For the rest of the readers, the
critical points of $f_n$ come in symmetric pairs.

\subsection{}
Instead of distances between critical points, we now have to talk about distances between
symmetric pairs of points on $\bS^2(n)$.

The values $f_n(p)$ and $f_n(-p)$ are equal up to the sign ($+$ if $f_n$ is even,
$-$ if $f_n$ is odd), while for $n^{1-C\e}$-separated pairs $(p_1, -p_1)$ and $(p_2, -p_2)$,
the random variables $f_n(p_1)$ and $f_n(p_2)$ are almost independent.

When applying Lemma~\ref{lemma:independence} to replace $g_n(p)$ by $\xi(p)$ we keep the relation between $\xi(p)$ and
$\xi (-p)$ the same as between $g_n(p)$ and $g_n(-p)$, i.e., they coincide up to a sign,
and make $\xi(p)$ and $\xi (p')$ independent for $p\ne \pm p'$.

\subsection{}
The rest of the argument goes as before with one simplification and two minor caveats.

\subsubsection{}
The simplification is that for the spherical harmonics ensemble, the blinking circles
cannot occur: by the classical Faber-Krahn inequality the area of any nodal domain
of a spherical harmonic on the sphere $\bS^2(n)$ cannot be less than a positive
numerical constant. So we need to treat only the Bogomolny-Schmit loops.

\subsubsection{}
Both caveats pertain to the proof of Lemma~\ref{lemma:variance_clusters} which estimates
from below the fluctuations in the number of random loops.
First of all, we note that since all steps of our construction were symmetric with respect to
the mapping $x\mapsto -x$, the results it produces are also symmetric. In particular the joints $J(p, \delta)$ and $J(-p, \delta)$ are symmetric,
and the set of connected components of
\[
Z(f_L)\setminus\bigcup_{p\in\Cr_{\tt S}(\alpha)} J(p, \delta)
\]
that touch the boundaries $\partial J(p, \delta)$ is also symmetric. Hence, the graph $G(V,E)$, to which Lemma~\ref{lemma:variance_clusters} was applied, is symmetric as well.

When defining marked cycles, we cannot choose a cycle $\ell$
passing through antipodal vertices, i.e., having common vertices with the symmetric
cycle $-\ell$. Fortunately, this does not happen often: there are at most $8$
faces $\mathfrak f$ such that $\mathfrak f$ and $-\mathfrak f$ have a common vertex $v$
(in which case, $-v$ will be also a common vertex of $\mathfrak f$ and $-\mathfrak f$).
This follows from the following lemma.

\begin{lemma}\label{Lemma:SphTopology}
Let $X, X' \subset \bS^2$ be two closed symmetric (with respect to the inversion
$x\mapsto -x$) non-empty symmetric sets. Then $X \cap X' \ne \emptyset$.
\end{lemma}

First, we conclude our argument, and then will prove Lemma~\ref{Lemma:SphTopology}.
Let $\mathfrak f$ and $-\mathfrak f$ have a pair of common vertices
$v$, $-v$ on their boundaries, then we can join $v$ and $-v$ by a path $\gamma$ with
$\gamma\setminus \{v, -v\}\subset \mathfrak f$, and by the path $-\gamma$ with
$-\gamma\setminus \{v, -v\}\subset -\mathfrak f$. Put
$X = \gamma\cup - \gamma$. This is a symmetric closed connected subset
of $\bS^2$. If $\mathfrak f'$ is another such face (different from $\pm \mathfrak f$), then we have another symmetric closed connected set
$X'$, and, by Lemma~\ref{Lemma:SphTopology},
$X \cap X' \ne\emptyset$.
Since $(\gamma\cup - \gamma)\setminus \{v, -v\}$ is contained inside
$\mathfrak f \cup - \mathfrak f$,
while $(\mathfrak f \cup - \mathfrak f) \cap (\mathfrak f' \cup - \mathfrak f') = \emptyset$,
we see that $v$ and $-v$ must be vertices on
$\partial \mathfrak f' \cup \partial(-\mathfrak f')$ as well. Recalling that each vertex on our graph has degree $4$, we see that there are at most
$8$ such ``bad faces'' $\mathfrak f$.

\medskip

\paragraph{\em Proof of Lemma~\ref{Lemma:SphTopology}} Assume that $X_1\cap X_2=\emptyset$. Without loss of generality, we assume that the set $X_2$ contains the North and the
South Poles. Since $X_1$, $X_2$ are compact, there exists $\e>0$ such that
$\operatorname{dist}(X_1, X_2)>6\e$.

Take any point $z\in X_1$. Then $-z\in X_1$ as well.
Since $X_1$ is connected, there exists a finite chain of points in $X_1$
$z=z_0, z_1, \ldots , z_n=-z$, with $d(z_j, z_{j-1})<\e$ ($j=1, \ldots , n$).
Connecting $z_{j-1}$ to $z_j$ by the shortest arc, we get a curve $\gamma_1$ joining
$z$ to $-z$ and staying in the $\e$-neighbourhood of $X_1$. Let
$\gamma = \gamma_1 \cup (-\gamma_1)$. Then $\gamma$ is a symmetric curve going from $z$ to $-z$ and back and staying in the $\e$-neighbourhood of $X_1$.

In a similar way, we can construct a curve from the North Pole to the South Pole staying in the $\e$-neighbourhood of
$X_2$. Let $\gamma_2$ be the piece of that  curve from the last intersection with the circle of radius $\e$ around the North Pole to the first intersection with the circle of radius $\e$ around the South Pole. Note that $\gamma_2$ and both these circles stay in the $\e$-neighbourhood of $X_2$ and, thereby, are disjoint with $\gamma$.

Now take the projection
\[
(x_1, x_2, x_3) \mapsto \Bigl( \frac{x_1}{\sqrt{x_1^2 + x_2^2}},
\frac{x_2}{\sqrt{x_1^2 + x_2^2}}, x_3 \Bigr)
\]
of $\bS^2\cap \{|x_3|\le \cos\e\}$ onto the cylinder $C=\{x_1^2+x_2^2=1, |x_3|\le \cos\e\}$.
Note that it preserves the symmetry with respect to the origin, so we get two disjoint
curves $\widetilde{\gamma}$ and $\widetilde{\gamma}_2$ on this cylinder such that
$ \widetilde{\gamma} = \widetilde{\gamma_1} \cup (-\widetilde{\gamma_1})$, where
$ \widetilde{\gamma_1} $ goes from some point $\widetilde z\in C$ to $-\widetilde z$ and stays at positive distance from the edge circles $\{x_3=\pm\cos\e\}$ of $C$, while $\widetilde{\gamma}_2$ joins those edge circles.

Consider the universal covering map ${\mathfrak p}\colon S\to C$, where
$S=\{(t, s)\in\bR^2\colon |s|\le\cos\e\}$ is a horizontal infinite strip
and ${\mathfrak p}((t,s))=(\cos2\pi t, \sin 2\pi t, s)$.
Note that, for $\ell\in\bZ$, ${\mathfrak p}((t+\tfrac12 +\ell, -s))
= - {\mathfrak p}((t, s))$. By the path lifting lemma,
$\widetilde{\gamma}_2$ is lifted to some curve $\Gamma_2$ on $S$
joining the top and the bottom boundary lines.
The curve $\widetilde{\gamma}_1$ is lifted to some curve
$\Gamma_1(\tau)$, $\tau\in [0, 1]$, joining $(t_0, s_0)$ with
$(t_0+\xi, -s_0)$ where $\xi\in\bZ+\tfrac12$. Then the curve
$\Gamma_1^*(\tau) = (t(\tau)+\xi, - s(\tau))$ extends $\Gamma_1$ and
projects to $-\widetilde{\gamma}_1$. This extension process can now be repeated
and done in both directions, so we get a curve $\Gamma$ in $S$ staying away
from the boundary and such that the first coordinate of $\Gamma$ goes from
$-\infty$ to $+\infty$ (if $\xi<0$, we reorient $\Gamma$). We still have
$\Gamma \cap \Gamma_2 = \emptyset$, so the increment of $\arg(w-w_2)$, as $w$ runs over $\Gamma$ and $w_2\in\Gamma_2$ stays fixed, should not depend on $w_2$.
However, this increment is $+\pi$ when $w_2$ is on the top boundary line of $S$ and
$-\pi$ when $w_2$ is on the bottom line. This contradiction proves the lemma.
\medskip{} \hfill $\Box$

\subsubsection{}
The second caveat is caused by the fact that
good cycles now come in symmetric pairs, and the cycles in each pair simultaneously
either merge other cycles or remain separate.
So our Bernoulli random variables $\eta_p$ are now valued
in $\{0, 2\}$ instead of $\{0, 1\}$.

\subsection*{Acknowledgement}
We are grateful to Dmitry Belyaev, Ron Peled, Evgenii Shustin, and Boris Tsirelson
for helpful discussions and suggestions.

\end{document}